\documentclass{amsart}
\usepackage{graphics,verbatim}	
\input epsf         

\newtheorem{proposition}{Proposition}[section]
\newtheorem{theorem}[proposition]{Theorem}

\newtheorem{lemma}[proposition]{Lemma}

\newtheorem{remark}[proposition]{Remark}
\newtheorem{example}[proposition]{Example}

\newtheorem{prop}[proposition]{Proposition}
\newtheorem{cor}[proposition]{Corollary}

\numberwithin{equation}{section}

\newcommand{\cov}{\mathrm{cov}}

\newcommand{\reals}{\mathbb R}

\newcommand{\cl}{\operatorname{cl}}
\newcommand{\nc}{\operatorname{nc}}
\newcommand{\Cat}{\operatorname{Cat}}

\newcommand{\covers}{{\,\,\,\cdot\!\!\!\! >\,\,}}

\newcommand{\set}[1]{{\left\lbrace #1 \right\rbrace}}

\newcommand{\Tge}{{T_{\ge -1}}}

\newcommand{\ep}{\varepsilon}
\newcommand{\br}[1]{\langle #1 \rangle}
\newcommand{\A}{{\mathcal A}}

\newcommand{\cm}{\parallel}

\newcommand{\toname}[1]{\stackrel{#1}{\longrightarrow}}

\newlength{\lsmash} % this sets how much the 2nd overline/underline in \up/\down is mashed down onto the first
\newlength{\lsmashalt} % for the alternate barring
\newlength{\myheight}
\newlength{\mydepth}

\newcommand{\up}[1]{
\settoheight{\myheight}{\ensuremath{\overline{#1}}}
\addtolength{\myheight}{-\lsmash}
\overline{
\protect \raisebox{0 pt}[\myheight][0 pt]{\ensuremath{\overline{#1}}}
}
}

\newcommand{\down}[1]{
\settodepth{\mydepth}{\ensuremath{\underline{#1}}}
\addtolength{\mydepth}{-\lsmash}
\underline{
\protect \raisebox{0 pt}[0 pt][\mydepth]{\ensuremath{\underline{#1}}}
}
}

\newcommand{\upwide}[1]{
\settoheight{\myheight}{\ensuremath{\overline{\,#1\,}}}
\addtolength{\myheight}{-\lsmash}
\overline{
\protect \raisebox{0 pt}[\myheight][0 pt]{\ensuremath{\overline{\,#1\,}}}
}
}

\newcommand{\downwide}[1]{
\settodepth{\mydepth}{\ensuremath{\underline{\,#1\,}}}
\addtolength{\mydepth}{-\lsmash}
\underline{
\protect \raisebox{0 pt}[0 pt][\mydepth]{\ensuremath{\underline{\,#1\,}}}
}
}

\begin{document}
\title[Coxeter-sortable elements]{Clusters, Coxeter-sortable elements and noncrossing partitions}

\author{Nathan Reading}
\address{
Mathematics Department\\
    University of Michigan\\
    Ann Arbor, MI 48109-1043\\
USA}
\thanks{The author was partially supported by NSF grant DMS-0202430.}
\email{nreading@umich.edu}
\urladdr{http://www.math.lsa.umich.edu/\textasciitilde nreading/}
%\keywords{}
%\date{\today}

\begin{abstract}
We introduce Coxeter-sortable elements of a Coxeter group~$W.$\
For finite $W,$\ we give bijective proofs that Coxeter-sortable elements are equi\-numerous with clusters and with noncrossing partitions.
We characterize Coxeter-sortable elements in terms of their inversion sets and, in the classical cases, in terms of permutations.
\end{abstract}

\maketitle

\setcounter{tocdepth}{1}
\tableofcontents

\section*{Introduction}
The famous Catalan numbers can be viewed as a special case of the \emph{\mbox{$W$-Catalan} number}, which counts various objects related to a finite Coxeter group~$W.$\
In many cases, the common numerology is the only known connection between different objects counted by the $W$-Catalan number.
One collection of objects counted by the $W$-Catalan number is the set of noncrossing partitions associated to~$W,$\ which play a role in low-dimensional topology, geometric group theory and non-commutative probability~\cite{McSurvey}.
Another is the set of maximal cones of the cluster fan.
The cluster fan is dual to the generalized associahedron~\cite{gaPoly,ga}, a polytope whose combinatorial structure underlies cluster algebras of finite type~\cite{CA2}.

This paper connects noncrossing partitions to associahedra via certain elements of~$W$ which we call \emph{Coxeter-sortable} elements or simply {\em sortable} elements.
For each Coxeter element $c$ of~$W,$\ there is a set of $c$-sortable elements, defined in the context of the combinatorics of reduced words.
We prove bijectively that sortable elements are equinumerous with clusters and with noncrossing partitions.
Sortable elements and the bijections are defined without reference to the classification of finite Coxeter groups, but the proof that these maps are indeed bijections refers to the classification.
The bijections are well-behaved with respect to the refined enumerations associated to the Narayana numbers and to positive clusters.

In the course of establishing the bijections, we characterize the sortable elements in terms of their inversion sets.
Loosely speaking, we ``orient'' each rank-two parabolic subgroup of~$W$ and require that the inversion set of the element be ``aligned'' with these orientations.
In particular, we obtain a characterization of the sortable elements in types $A_n$, $B_n$ and $D_n$ as permutations.

Because sortable elements are defined in terms of reduced words, it is natural to partially order them as an induced subposet of the weak order.
Indeed, the definition of sortable elements arose from the study of certain lattice quotients of the weak order called Cambrian lattices.
%sinceV1:  The following sentence altered:
In the sequel~\cite{sort_camb} to this paper, we show that sortable elements are a combinatorial model for the Cambrian lattices.

Recently, Brady and Watt~\cite{BWlattice} observed that the cluster fan arises naturally in the context of noncrossing partitions.
Their work and the present work constitute two seemingly different approaches to connecting noncrossing partitions to clusters.
The relationship between these approaches is not yet understood.

The term ``sortable'' has reference to the special case where~$W$ is the symmetric group.
For one choice of~$c$, the definition of $c$-sortable elements of the symmetric group is as follows:
Perform a bubble sort on a permutation~$\pi$ 
by repeatedly passing from left to right in the permutation and, whenever two adjacent elements are out of order, transposing them.
For each pass, record the set of positions at which transpositions were performed.
Then~$\pi$ is~$c$-sortable if this sequence of sets is weakly decreasing in the containment order.
There is also a choice of $c$ (see Example~\ref{Ansort}) such that the $c$-sortable elements are exactly the $231$-avoiding or stack-sortable permutations \cite[Exercise 2.2.1.4--5]{Knuth}.

\section{Coxeter groups}
\label{cox sec}
In this section we establish conventions and present some basic lemmas about Coxeter groups.
We assume standard notions from the theory of Coxeter groups, as found for example in~\cite{Bourbaki,Humphreys}.

\subsection*{Coxeter arrangements and parabolic subgroups}
Throughout this paper,~$W$ denotes a finite Coxeter group with simple generators~$S$ and reflections $T$.
Some definitions apply also to the case of infinite $W,$\ but we confine the treatment of the infinite case to a series of remarks (Remarks~\ref{inf sort}, \ref{inf nc} and~\ref{inf cl}.).

For every definition and result in this paper, without exception, the relationship between the case of irreducible~$W$ and the case of general~$W$ is completely straightforward.
In particular, for those results which are proved using the classification of finite Coxeter groups, we pass to the irreducible case without comment.

The term ``word'' always means ``word in the alphabet~$S$.''
Later, we consider words in the alphabet $T$ which, to avoid confusion, are called ``$T$-words.''
An {\em inversion} of $w\in W$ is a reflection $t$ such that $l(tw)<l(w)$.
The set of inversions of~$w$ is denoted $I(w)$.
The longest element of~$W$ is denoted~$w_0$.
A {\em descent} of~$w$ is a simple reflection~$s$ such that $l(ws)<l(w)$.
A {\em cover reflection} of $w\in W$ is a reflection $t$ such that $tw=ws$ for some descent~$s$ of~$w$.
The term ``cover reflection'' refers to the (right) weak order.
This is the partial order on $W$ whose cover relations are the relations of the form $w\covers ws$ for a descent~$s$ of~$w$, or equivalently, $w\covers tw$ for a cover reflection $t$ of~$w$.

For a word $a=a_1a_2\cdots a_k$, define the {\em support} of $a$ to be $\set{a_1,a_2,\ldots,a_k}$.
For $w\in W$, the support of~$w$ is the support of any reduced word for~$w$, and~$w$ is said to have {\em full support} if its support is~$S$.

A finite Coxeter group~$W$ has a representation as a group generated by orthogonal reflections in $\reals^n$ for $n=\mbox{rank}(W)$.
We assume that some reflection representation has been fixed and identify the elements of~$W$ with orthogonal transformations
via this representation.
The set of reflecting hyperplanes for the reflections of~$W$ is called the {\em Coxeter arrangement} $\A$ associated to~$W.$\
The {\em regions} of $\A$ are the closures of the connected components of $\reals^n-\cup\A$.
Given the group~$W$ represented as a reflection group, one chooses a set~$S$ of simple generators by choosing some region $R$ and 
taking~$S$ to be the set of reflections in the facet hyperplanes of $R$.
The regions of $\A$ are in one-to-one correspondence with the elements of~$W,$\ with an element~$w$ corresponding to the region~$w(R)$.
The inversion set of~$w$ corresponds to the set of reflecting hyperplanes which separate~$w(R)$ from $R$.
The cover reflections of~$w$ correspond to the inversions of~$w$ which are also facet-defining hyperplanes of~$w(R)$.
For each $s\in S$ there is a unique extreme ray $\rho(s)$ of $R$ which is not contained in the facet hyperplane corresponding to~$s$.

Given a subset $J\subseteq S$, let~$W_J$ denote the subgroup of~$W$ generated by $J$.
This is a {\em standard parabolic subgroup} of~$W.$\
Often we consider $J=\br{s}$ where $\br{s}$ stands for $S-\set{s}$.
A reflection in~$W$ is in~$W_J$ if and only if the corresponding hyperplane contains $\rho(s)$ for every $s\not\in J$.
An element $w\in W$ is in~$W_J$ if and only if~$w$ fixes $\rho(s)$ for every $s\not\in J$.
Equivalently, $w\in W_J$ if and only if the region corresponding to~$w$ contains $\rho(s)$ for every $s\not\in J$.
In particular we have the following.

\begin{lemma}
\label{cover para}
For any $J\subseteq S$, an element $w\in W$ is in~$W_J$ if and only if every cover reflection of~$w$ is in~$W_J$.
\end{lemma}

A {\em parabolic subgroup} is any subgroup~$W'$ conjugate to a standard parabolic subgroup.
%sinceV1:changed following sentences
Our terminology for parabolic subgroups is consistent with~\cite{Bourbaki} and~\cite{Humphreys}, but the term ``parabolic subgroup'' is often used for what we are here calling a ``standard parabolic subgroup.''
We define {\em canonical generators} for a (not necessarily standard) parabolic subgroup~$W'$ as follows:
The set~$S$ of simple generators for~$W$ arises from the choice of a region $R$ as explained above.
The canonical generators for~$W'$ are the simple generators for~$W'$ obtained by choosing a region $R'$ in the Coxeter arrangement for~$W'$ 
such that $R\subseteq R'$.
The following two lemmas are immediate.

\begin{lemma}
\label{s in}
Let~$W'$ be a parabolic subgroup with canonical generators $S'$.
If $s\in S\cap W'$ then $s\in S'$.
\end{lemma}

\begin{lemma}
\label{canon lem}
For any $w\in W$, the cover reflections of~$w$ generate a parabolic subgroup of~$W$ and are the canonical generators of that parabolic subgroup.
\end{lemma}

Let $W'$, $S'$, $R'$ and $R$ be as above, with $R\subseteq R'$.
For $s\in S$, the cone $s(R')$ is a region in the Coxeter arrangement for $sW's$, and contains $R$ if and only if $s\not\in S'$.
Thus the following statement holds.

\begin{lemma}
\label{canon conj}
Let~$W'$ be a parabolic subgroup with canonical generators $S'$ and let $s\in S$.
If $s\not\in S'$ then the parabolic subgroup $sW's$ has canonical generators $sS's=\set{ss's:s'\in S'}$.
\end{lemma}
The following lemma is a special case of \cite[Lemma~6.6]{congruence}.

\begin{lemma}
\label{int}
For any standard parabolic subgroup~$W_J$ and any parabolic subgroup~$W'$ of rank two, if the intersection $W'\cap W_J$ is
nonempty then either it is a single canonical generator of~$W'$ or it is all of~$W'.$\
\end{lemma}

The following is part of \cite[Lemma~2.11]{Dyer}.

\begin{lemma}
\label{dyer lem}
Let~$W'$ be a rank-two parabolic subgroup with canonical generators $t_1$ and $t_2$ and let $w\in W$.
Then $I(w)\cap W'$ is either an initial segment or a final segment of the sequence 
\[t_1,\,t_1t_2t_1,\,\ldots,\,t_2t_1t_2,\,t_2.\]
\end{lemma}

\subsection*{Coxeter elements}
For the rest of the paper,~$c$ denotes a {\em Coxeter element}, that is, an element of~$W$ with a reduced word which is a permutation of~$S$.
An {\em orientation} of the Coxeter diagram for~$W$ is obtained by replacing each edge of the diagram by a single directed edge, connecting the same pair of vertices in either direction.
Orientations of the Coxeter diagram correspond to Coxeter elements (cf.~\cite{shi}) as follows:
Given a Coxeter element~$c$, any two reduced words for~$c$ are related by commutations of simple generators.
An edge $s$---$t$ in the diagram represents a pair of noncommuting simple generators, and the edge is oriented $s\longrightarrow t$ if and only if~$s$ precedes $t$ in every reduced word for~$c$.

A simple generator $s\in S$ is {\em initial} in (or is an {\em initial letter} of) a Coxeter element~$c$ if it is the first letter of some reduced word for~$c$.
Similarly~$s$ is {\em final} in~$c$ if it is the last letter  of some reduced word for~$c$.
If~$s$ is initial in~$c$ then $scs$ is a Coxeter element and~$s$ is final in $scs$.
%Since any two reduced words for~$c$ are related by commutations of simple generators, any two initial letters of~$c$ commute.
%Similarly any two final letters of~$c$ commute.
For any Coxeter element~$c$ of~$W$ and any $J\subseteq S$, the {\em restriction} of~$c$ to~$W_J$ is the Coxeter element for~$W_J$ obtained
by deleting the letters $S-J$ from any reduced word for~$c$.

The following is a special case of \cite[Theorem~1.2(1)]{BGP}.

\begin{lemma}
\label{swap sources}
Given any two Coxeter elements~$c$ and $c'$ of~$W,$\ there is a sequence $c=c_0,c_1,\ldots,c_k=c'$ such that 
for each $i\in[k]$ there is an initial letter $s_{i-1}$ of $c_{i-1}$ such that $c_i=s_{i-1}c_{i-1}s_{i-1}$.
\end{lemma}

The Coxeter diagram of any finite Coxeter group is a labeled bipartite graph.
Let $S_+$ and $S_-$ be the parts in any bipartition of the diagram.
A \emph{bipartite} Coxeter element is $c_-c_+$, where 
\[c_+:=\prod_{s\in S_+}s\qquad\mbox{and}\qquad c_-:=\prod_{s\in S_-}s.\]
Since elements of $S_\ep$ commute pairwise for $\ep\in\set{+,-}$, these products are well-defined.
The following lemma is part of \cite[Lemma~1.3.4]{Bessis}.

\begin{lemma}
\label{c orbit}
For any $t\in T$, the orbit of $t$ under conjugation by the dihedral group $\br{c_-,c_+}$ contains an element of~$S$.
\end{lemma}

\begin{remark}\rm
Throughout the paper, various constructions involve a choice of Coxeter element~$c$.  
The choice of~$c$ is sometimes indicated by a subscript.
We also use these explicit references to~$c$ as a way of specifying a standard parabolic subgroup.
For example, if~$s$ is initial in~$c$ then $sc$ is a Coxeter element in the standard parabolic subgroup~$W_{\br{s}}$.
Thus an explicit reference to $sc$ in a construction implies in particular that the construction takes place in~$W_{\br{s}}$.
\end{remark}

The following lemma provides an inductive structure for the connection between sortable elements and noncrossing partitions in Section~\ref{nc sec}.
Let~$W'$ be a nontrivial parabolic subgroup with canonical generators $S'$ and let~$c$ be a Coxeter element.
For any sequence $s_0,s_1,\ldots,s_k$ of simple generators, define $W'_0=W'$, $c_0=c$, $W'_i=s_{i-1}W'_{i-1}s_{i-1}$ and 
$c_i=s_{i-1}c_{i-1}s_{i-1}$.
Let $S'_i$ be the set of canonical generators for~$W'_i$.

\begin{lemma}
\label{nu lemma}
For any parabolic subgroup~$W'$ there exists a sequence $s_0,s_1,\ldots,s_k$ with $s_i$ initial in $c_i$ for $i=0,1,\ldots,k$ and $S'_i=s_{i-1}S'_{i-1}s_{i-1}$ 
for each $i\in[k]$, such that $s_k\in S'_k$.
\end{lemma}
\begin{proof}
By Lemma~\ref{swap sources} there is a sequence $s_0,s_1,\ldots,s_{j-1}$ of simple generators such that $s_i$ is initial in $c_i$ for $i=0,1,\ldots,j-1$
and such that $c_j$ is a bipartite Coxeter element with respect to some bipartition $S=S_+\cup S_-$.
If there is some $i$ with $0\le i\le j-1$ such that $s_i\in S'_i$ then we are done.
Otherwise, let $\Sigma_0=S'_j$, $\Sigma_1=(c_j)_-\Sigma_0(c_j)_-$, $\Sigma_2=(c_j)_+\Sigma_1(c_j)_+$, $\Sigma_3=(c_j)_-\Sigma_2(c_j)_-$, and so forth.
Lemma~\ref{c orbit} says that for some $m\ge 0$, $\Sigma_m$ contains a simple generator~$s$. 
Since $(c_j)_\ep$ fixes each generator in $S_\ep$, we furthermore choose the smallest $m\ge 0$ such that both $\Sigma_m$ and $\Sigma_{m+1}$ contain~$s$.
For some particular reduced words for $(c_j)_+$ and $(c_j)_-$, we append the letters of $(c_j)_+$ to the sequence $s_0,s_1,\ldots,s_{j-1}$ followed by the
letters of $(c_j)_-$, then the letters of $(c_j)_+$, and so forth until we have appended the letters of $(c_j)_-(c_j)_+(c_j)_-\cdots(c_j)_{(-1)^m}$  ($m$ factors).
Calling the sequence thus obtained  $s_0,s_1,\ldots,s_{k-1}$ and setting $s_k=s$, we have a sequence with the desired properties.
\end{proof}

Define $\nu_c(W')$ to be the smallest $k$ such that there exists a sequence $s_0,s_1,\ldots,s_k$ satisfying the conclusions of 
Lemma~\ref{nu lemma}.
In particular, $\nu_c(W')=0$ if and only~$W'$ contains some initial letter $s$ of $c$.

\subsection*{The classical types}
We now review the usual combinatorial realizations of the Coxeter groups of types A, B and D\@.
We also introduce notation which is used in Theorem~\ref{ABDsort} to characterize $c$-sortable elements of these types.

We begin with the case $W=A_n$, which is realized combinatorially as the symmetric group $S_{n+1}$.
Permutations $\pi\in S_{n+1}$ are written in {\em one-line notation} as $\pi_1\pi_2\cdots\pi_{n+1}$ with $\pi_i=\pi(i)$.
The simple generators of $S_{n+1}$ are $s_i=(i\,\,\,\,i\!+\!1)$ for $i\in[n]$.
Reflections are transpositions $(i\,\,j)$ for $1\le i<j\le n+1$.
A descent of~$\pi$ is $s_i$ such that $\pi_i>\pi_{i+1}$, and the corresponding cover reflection is $(\pi_{i+1}\,\,\pi_i)$.
A~simple generator $s_i$ has $l(s_i\pi)<l(\pi)$ if and only if $i+1$ precedes $i$ in the one-line notation for~$\pi$.
A permutation~$\pi$ is in the standard parabolic subgroup~$W_{\br{s_i}}$ if and only if every element of $[i]$ precedes every element of $[i+1,n+1]$ in the one-line notation for~$\pi$.
A reflection $(j\,\,k)$ with $j<k$ is in~$W_{\br{s_i}}$ if and only if either $k\le i$ or $j\ge i+1$.

A {\em barring} of a set $U$ of integers is a partition of that set into two sets $\up{U}$ and $\down{U}$.
Elements of $\up{U}$ are {\em upper-barred} integers denoted $\upwide{i}$ and {\em lower-barred} integers are elements of $\down{U}$, denoted $\downwide{i}$.

\begin{figure}[ht]
\centerline{
\scalebox{.8}{
\epsfbox{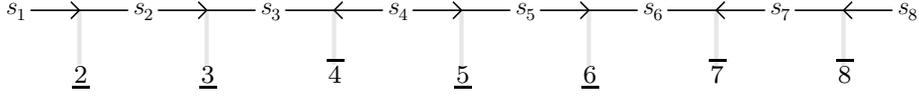}
}
}
\caption{Orientation and barring in $S_9$}
\label{barA}
\end{figure}

Recall that orientations of the Coxeter diagram correspond to Coxeter elements.
The Coxeter diagram for $S_{n+1}$ has unlabeled edges connecting $s_i$ to $s_{i+1}$ for $i\in[n-1]$. 
We encode orientations of the Coxeter diagram for $S_{n+1}$ as barrings of $[2,n]$ by directing $s_i\to s_{i-1}$ for every 
$\upwide{i}\in[2,n]$ and $s_{i-1}\to s_i$ for every $\downwide{i}\in [2,n]$, as illustrated in Figure~\ref{barA} for 
$c=s_8s_7s_4s_1s_2s_3s_5s_6$ in $S_9$.
Given a choice of Coxeter element, the corresponding barring is assumed.

We now consider the case $W=B_n$, employing the usual combinatorial realization of $B_n$ as the group of {\em signed permutations} 
of $\pm[n]=\set{\pm 1,\pm 2,\ldots,\pm n}$.
These are permutations~$\pi$ of $\pm[n]$ with $\pi(-i)=-\pi(i)$ for all $i\in[n]$.
We write these permutations in {\em long one-line notation} as $\pi_{-n}\pi_{-n+1}\cdots\pi_{-1}\pi_1\pi_2\cdots\pi_n$ with
\mbox{$\pi_i=\pi(i)$}.
The simple generators of $B_n$ are $s_0=(-1\,\,\,1)$ and $s_i=(-i\!-\!1\,\,\,-\!i)(i\,\,\,\,i\!+\!1)$ for $i\in[n-1]$.
The reflections in $B_n$ are transpositions $(-i\,\,i)$ for $i\in[n]$ and pairs of transpositions $(i\,\,j)(-j\,\,-i)$ for $i\neq -j$.
We refer to the latter type of reflection by specifying only one of the two transpositions.
A descent of a signed permutation~$\pi$ is $s_i$ such that $\pi_{i+1}$ is less than the preceding entry in the long one-line notation for~$\pi$.
The corresponding cover reflection transposes $\pi_{i+1}$ with the preceding entry.
For $i>0$, a signed permutation~$\pi$ has $l(s_i\pi)<l(\pi)$ if and only if $i+1$ precedes $i$ in the one-line notation for~$\pi$.
The signed permutation has $l(s_0\pi)<l(\pi)$ if and only if $1$ precedes $-1$.
A signed  permutation~$\pi$ is in the standard parabolic subgroup~$W_{\br{s_i}}$ if and only if every element of $[-n,i]-\set{0}$ precedes every element of $[i+1,n]$ in the long one-line notation for~$\pi$.
A reflection $(j\,\,k)$ with $|j|\le k$ is in~$W_{\br{s_i}}$ if and only if either $k\le i$ or $j\ge i+1$.

\begin{figure}[ht]
\centerline{
\scalebox{.8}{
\epsfbox{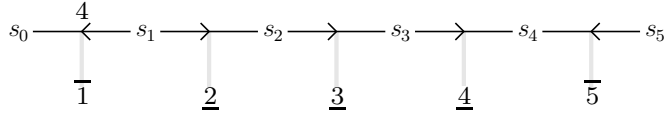}
}
}
\caption{Orientation and barring in $B_6$}
\label{barB}
\end{figure}

The Coxeter diagram for $B_n$ has unlabeled edges connecting $s_i$ to $s_{i+1}$ for $i\in[n-2]$ and an edge labeled $4$ connecting $s_0$ to $s_1$. 
Orientations of the Coxeter diagram for $B_n$ are encoded by barrings of $[n-1]$ by directing $s_i\to s_{i-1}$ for every 
$\upwide{i}\in[n-1]$ and $s_{i-1}\to s_i$ for every $\downwide{i}\in [n-1]$, as illustrated in Figure~\ref{barB} for $c=s_5s_1s_0s_2s_3s_4$ 
in $B_6$.
The barring of $[n-1]$ is extended to a barring of $\pm[n-1]$ by specifying that the barring of $-i$ is opposite the barring of $i$.

We now consider the Coxeter group $D_n$ with the usual combinatorial realization as the group of {\em even-signed permutations}.
These are signed permutations~$\pi$ of $\pm[n]$ such that the set $\set{i\in[n]:\pi_i<0}$ has an even number of elements.
The simple generators are $s_i=(i\,\,\,i+1)(-i\!-\!1\,\,\,-i)$ for $i\in[n-1]$, and $s_0=(-2\,\,\,1)(-1\,\,\,2)$.
The reflections are pairs of transpositions $(i\,\,\,j)(-j\,\,\,-i)$ for $i\neq -j$.
We refer to reflections by naming one of the two transpositions.
For $i>0$, $s_i$ is a descent of $\pi\in D_n$ if and only if $\pi_i>\pi_{i+1}$.
The corresponding cover reflection is $(\pi_{i+1}\,\,\pi_i)$.
The generator $s_0$ is a descent if and only if $\pi_{-1}>\pi_2$.
The corresponding cover reflection is $(\pi_2\,\,\pi_{-1})$.
For $i>0$, $l(s_i\pi)<l(\pi)$ if and only if $i+1$ precedes $i$ in the long one-line notation for~$\pi$ and $l(s_0\pi)<l(\pi)$ is and only if $2$ precedes $-1$.
%For $i\neq 1$, an even-signed permutation~$\pi$ has $\pi\in W_{\br{s_i}}$ if and only if every element of $[-n,i]-\set{0}$ precedes every element of $[i+1,n]$ in the long one-line notation for~$\pi$.
%For $i\neq 1$, a reflection $(j\,\,k)$ with $|j|<k$ is in~$W_{\br{s_i}}$ if and only if either $k\le i$ or $j\ge i+1$.
%An even-signed permutation~$\pi$ is in $W_{\br{s_1}}$ if and only if it permutes the set $\set{-1}\cup[2,n]$.

\begin{figure}[ht]
\centerline{
\scalebox{.8}{
\epsfbox{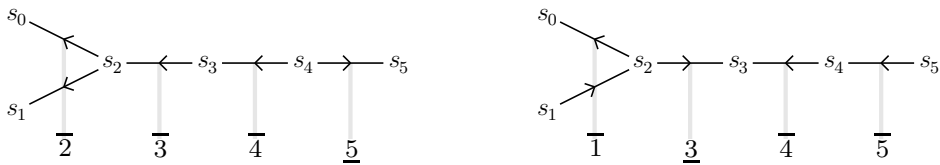}
}
}
\caption{Orientation and barring in $D_6$}
\label{barD}
\end{figure}

The Coxeter diagram for $D_n$ has unlabeled edges connecting $s_i$ to $s_{i+1}$ for $i\in[n-2]$
and an unlabeled edge connecting $s_0$ to $s_2$.
The Coxeter elements of $D_n$ are divided into two classes.
The  {\em symmetric Coxeter elements} are those which are preserved under the diagram automorphism which exchanges $s_0$ and $s_1$ 
and fixes all other simple generators.
{\em Asymmetric Coxeter elements} are not preserved by this automorphism.
A symmetric Coxeter element corresponds to a barring of $[2,n-1]$ by directing $s_i\to s_{i-1}$ for every $\upwide{i}\in[3,n-1]$ and 
$s_{i-1}\to s_i$ for every $\downwide{i}\in [3,n-1]$.
Additionally, if $2$ is lower-barred then direct $s_0\to s_2$ and $s_1\to s_2$ and if $2$ is upper-barred then direct $s_2\to s_0$ and $s_2\to s_1$.
An asymmetric Coxeter element corresponds to a barring of $\set{1}\cup[3,n-1]$.
The barring of $[3,n-1]$ is used to direct edges exactly as in the symmetric case.
If $1$ is lower-barred then direct $s_0\to s_2$ and $s_2\to s_1$ and if $1$ is upper-barred then direct $s_2\to s_0$ and $s_1\to s_2$.
The correspondence between orientations and barrings is illustrated in Figure~\ref{barD}.
The barring of $[2,n-1]$ (respectively $\set{1}\cup[3,n-1]$) is extended to a barring of $\pm[2,n-1]$ (respectively 
$\set{\pm1}\cup\pm[3,n-1]$) by requiring that $-i$ and $i$ have opposite barrings.

\begin{remark}\rm
\label{Dn res}
We introduce restriction of even-signed permutations.
Let $W$ be of type $D_n$ and let $X=\set{x_1,\ldots,x_{n'}}$ where the $x_i$ satisfy $0<x_1<\cdots<x_{n'}\le n$.
Given $\pi\in W$, consider the signed permutation obtained by restricting~$\pi$ to $\pm X$.
If this is not even-signed, then it can be made even-signed by transposing the two centermost entries in its long one-line notation.
We obtain an element $\pi'\in W'$, where $W'$ is of type $D_{n'}$ and is realized as signed permutations of $\pm X$ with simple generators $s'_i=(-a_{i+1}\,\,\,-{a_i})(a_i\,\,\,a_{i+1})$ for $1\le i\le n'-1$ and $s'_0=(-a_1\,\,\,a_2)(-a_2\,\,\,a_1)$.
If $\set{1,2}\subseteq X$ then the barring defining a Coxeter element $c$ of $W$ restricts to a barring of certain elements of $\pm X$ so as to uniquely define a Coxeter element $c'$ of $W'$.
This restriction operation allows us to reduce certain proofs (Lemmas~\ref{align sc} and~\ref{nc lemma}) in type D to the case of low rank, and check the low-rank examples by computer.
In effect, this replaces a many-case prose proof by a computer check of the cases.
\end{remark}

\section{Coxeter-sortable elements}
\label{sort sec}
In this section we define Coxeter-sortable elements and present two easy lemmas which play a key role throughout the paper.
Fix a reduced word $a_1a_2\cdots a_n$ for a Coxeter element~$c$.
Write a half-infinite word 
\[c^\infty=a_1a_2\cdots a_na_1a_2\cdots a_na_1a_2\cdots a_n\ldots\]
The {\em $c$-sorting word} for $w\in W$ is the lexicographically first (as a sequence of positions in $c^\infty$) subword of $c^\infty$ 
which is a reduced word for~$w$.
The $c$-sorting word can be interpreted as a sequence of subsets of~$S$ by rewriting
\[c^\infty=a_1a_2\cdots a_n|a_1a_2\cdots a_n|a_1a_2\cdots a_n|\ldots\]
where the symbol ``$|$'' is called a {\em divider}.
The subsets in the sequence are the sets of letters of the $c$-sorting word which occur between adjacent dividers.
This sequence contains a finite number of non-empty subsets, and furthermore if any subset in the sequence is empty, then every 
later subset is also empty.
For clarity in examples, we often retain the dividers when writing $c$-sorting words for $c$-sortable elements.

An element $w\in W$ is {\em $c$-sortable} if its $c$-sorting word defines a sequence of subsets which is decreasing under inclusion.
This definition of $c$-sortable elements requires a choice of reduced word for~$c$.
However, for a given~$w$, the $c$-sorting words for~$w$ arising from different reduced words for~$c$ are related by commutations of letters, with no commutations across dividers.
In particular, the set of $c$-sortable elements does not depend on the choice of reduced word for~$c$.

\begin{remark}\rm
\label{tree rem}
The $c$-sortable elements have a natural search-tree structure rooted at the identity element.
The edges are pairs $v,w$ of $c$-sortable elements such that the $c$-sorting word for $v$ is obtained from the 
$c$-sorting word for~$w$ by deleting the rightmost letter.
This makes possible an efficient traversal of the set of \mbox{$c$-sortable} elements of~$W$ which, although it does not explicitly appear in what 
follows, allows various properties of $c$-sortable elements (in particular, Lemmas~\ref{nc lemma} and \ref{nc lemma 2}) to be checked computationally in low ranks.
Also, in light of the bijections of Theorems~\ref{nc} and~\ref{cl}, an efficient traversal of the $c$-sortable elements leads to an efficient traversal of noncrossing partitions or of clusters.
\end{remark}

\begin{example}\rm
\label{B2sort}
Consider $W=B_2$ with $c=s_0s_1$.
The $c$-sortable elements are $1$, $s_0$, $s_0s_1$, $s_0s_1|s_0$, $s_0s_1|s_0s_1$ and $s_1$.
The elements $s_1|s_0$ and $s_1|s_0s_1$ are not $c$-sortable.
\end{example}

\begin{example}\rm
\label{Ansort}
As a special case of Theorem~\ref{ABDsort}, when $W$ is $S_{n+1}$ and $c$ is the permutation \mbox{$(n\,\,\,\,n\!+\!1)\cdots(2\,\,3)(1\,\,2)$}, 
the $c$-sortable elements are exactly the $231$-avoiding or stack-sortable permutations defined in \cite[Exercise 2.2.1.4--5]{Knuth}.
\end{example}

The next two lemmas are immediate from the definition of $c$-sortable elements.
Together with the fact that $1$ is $c$-sortable for any $c$, they completely characterize $c$-sortability.

\begin{lemma}
\label{sc}
Let~$s$ be an initial letter of~$c$ and let $w\in W$ have $l(sw)>l(w)$.
Then~$w$ is $c$-sortable if and only if it is an $sc$-sortable element of~$W_{\br{s}}$.
\end{lemma}

\begin{lemma}
\label{scs}
Let~$s$ be an initial letter of~$c$ and let $w\in W$ have $l(sw)<l(w)$.
Then~$w$ is $c$-sortable if and only if $sw$ is $scs$-sortable.
\end{lemma}

\begin{remark}\rm
In the dictionary between orientations of Coxeter diagrams (i.e.\ quivers) and Coxeter elements, the operation
of replacing~$c$ by $scs$ corresponds to the operation on quivers which changes a source into a sink by reversing
all arrows from the source.
This operation was used in~\cite{MRZ} in generalizing the clusters of~\cite{ga} to $\Gamma$-clusters, where $\Gamma$ is a quiver of finite type.
We thank Andrei Zelevinsky for pointing out the usefulness of this operation, which plays a key role throughout the paper.
\end{remark}

\begin{remark}\rm
\label{inf sort}
The definition of $c$-sortable elements is equally valid for infinite $W$.
Lemmas~\ref{sc} and~\ref{scs} are valid and characterize $c$-sortability in the infinite case as well.
However, we remind the reader that for all stated results in this paper, $W$ is assumed to be finite.
\end{remark}

\section{Orienting a Coxeter group}
\label{or sec}
In this section we define the $c$-orientation of a finite Coxeter group~$W$ in preparation for the characterization of sortable elements in terms of their inversion sets in Theorem~\ref{align}.
At first glance the use of the term ``orientation'' here may appear to conflict with the use of that term in the context of orientations of Coxeter diagrams.
However, the $c$-orientation of~$W$ can be viewed as an extension of the diagram orientation (see Remark~\ref{exten}).

We remind the reader that throughout the paper, parabolic subgroups are not assumed to be necessarily standard.
Two reflections $t_1\neq t_2$ in~$W$ are {\em adjacent} if there is some region $R$ of the associated Coxeter arrangement such that
the reflecting hyperplanes of $t_1$ and $t_2$ define facets of $R$.
Equivalently, if~$W'$ is the rank-two parabolic subgroup containing $t_1$ and $t_2$, the two reflections are adjacent if there is some region $R'$ of the Coxeter arrangement for~$W'$ such that $t_1$ and $t_2$ define facets of $R'$.
A rank-two parabolic subgroup is irreducible if and only if it has more than two reflections, if and only if its two canonical generators do not commute.

An {\em orientation} of~$W$ is a directed graph with vertex set $T$ such that $t_1\to t_2$ implies that $t_1$ and $t_2$ are adjacent and noncommuting, with the property that for any irreducible rank-two parabolic subgroup~$W',$\ the directed graph induced on the reflections of~$W'$ is a directed cycle.
There are thus two distinct orientations of each irreducible rank-two parabolic subgroup and choosing an orientation of~$W$ is equivalent to choosing, independently, an orientation of each irreducible rank-two parabolic subgroup~$W'$.
Furthermore, the orientation of each irreducible rank-two parabolic is determined by any one of its directed edges.
In particular, one can specify orientations by specifying the directed edge between the canonical generators of each irreducible rank-two parabolic.

\begin{prop}
\label{family}
There is a unique family of orientations of~$W$ indexed by Coxeter elements and satisfying the following requirements (with directed arrows in the orientation indexed by~$c$ written $t_1\toname{c}t_2$).
For any irreducible rank-two parabolic subgroup~$W'$ and any initial letter~$s$ of~$c$:
\begin{enumerate}
\item[(i) ]If the canonical generators of~$W'$ are~$s$ and $t$ then $s\toname{c} t$, and 
\item[(ii) ]If $t_1\toname{c}t_2$ then $st_1s\toname{scs}st_2s$.
\end{enumerate}
The orientation indexed by~$c$ is called the {\em $c$-orientation} of~$W.$\
\end{prop}
\begin{proof}
In light of Lemma~\ref{nu lemma}, there is at most one family of orientations of~$W$ satisfying (i) and (ii).

We now introduce an explicit family of orientations.
Let $a_1\cdots a_m$ be a $c$-sorting word for~$w_0$.
Each reflection is $a_1a_2\cdots a_ja_{j-1}\cdots a_1$ for a unique $j\in[m]$, and we totally order $T$ by increasing $j$.
This is a {\em reflection order} in the sense of \mbox{\cite[Section~2]{Dyer}}.
This order induces an orientation of each irreducible rank-two parabolic subgroup~$W'$ as follows:
Let~$W'$ have canonical generators $t_1$ and $t_2$ such that $t_1$ precedes $t_2$ in the reflection order and declare that $t_1\toname{c}t_2$.
This orientation is independent of the choice of reduced word for~$c$ because changing the reduced word for~$c$ only changes the reflection order by transposition of commuting reflections.

Property (i) of this family of orientations is immediate for any~$W.$\
Property (ii) is trivial in rank two and is verified computationally for all of the exceptional types and for the classical types of low rank.
%Note that we need to do this for all irreducible types which are standard parabolic subgroups of exceptional groups. 
%To do this in sort_maple, run:
%crys:=true: for R in [A3,B3,A4,D4,B4,F4,A5,D5,B5,A6,D6,B6,E6,A7,D7,B7,E7,E8] do check_family_all(R) od;
%crys:=false: for R in [H3,H4] do check_family_all(R) od;
Presumably (ii) holds in all ranks for the classical types, but we have no proof.
Instead we prove Proposition~\ref{family} in the classical types by realizing the $c$-orientation in terms of planar diagrams.
The construction occupies most of the remainder of this section.
\end{proof}

In our combinatorial descriptions of the $c$-orientations of~$W$ in the cases $W=A_n$, $B_n$ or $D_n$, we continue the representation of Coxeter elements by barrings, as explained near the end of Section~\ref{cox sec}.

\begin{figure}[ht]
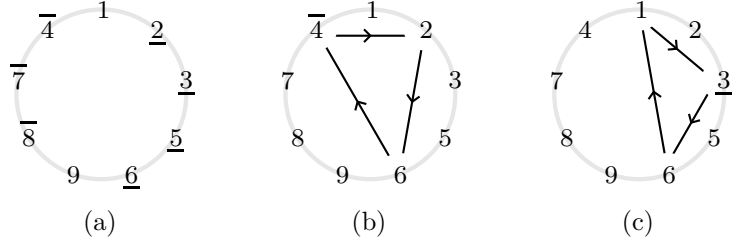

\centerline{
\begin{tabular}{ccccccc}
\scalebox{.8}{
\epsfbox{cycleAa.ps}
}
&&&
\scalebox{.8}{
\epsfbox{cycleAb.ps}
}
&&&
\scalebox{.8}{
\epsfbox{cycleAc.ps}
}
\\[4 pt]
(a)&&&(b)&&&(c)
\end{tabular}
}
\caption{$c$-Orientation of $S_9$}
\label{cycleA}
\end{figure}

A Coxeter element~$c$ of $S_{n+1}$ is realized as a permutation with a single cycle of the form
\[\left(1\,\,d_1\,\,d_2\,\,\cdots\,\,d_l\,\,(n+1)\,\,u_k\,\,u_{k-1}\,\,\cdots\,\,u_1\right),\]
where $d_1,d_2,\ldots d_l$ are the lower-barred elements of $[2,n]$ in increasing order and $u_1,u_2,\ldots,u_k$ are the upper-barred 
elements of $[2,n]$ in increasing order.
To model the $c$-orientation of~$W,$\ one first places the symbols in $[n+1]$ on a circle clockwise in the order given
by this cycle, as shown in Figure~\ref{cycleA}(a) for $c=s_8s_7s_4s_1s_2s_3s_5s_6$ in $S_9$ (cf. Figure~\ref{barA}).

Two noncommuting transpositions $t_1$ and $t_2$ in $S_{n+1}$ are the canonical generators of the rank-two parabolic subgroup~$W'$ which 
contains them if and only if one is $(i\,\,j)$ and the other is $(j\,\,k)$ with $i<j<k$.
These reflections, and the third reflection $(i\,\,k)$ in~$W'$ are naturally associated to the sides of an inscribed triangle in the cycle~$c$.
The $c$-orientation of~$W'$ is defined to be the cyclic ordering on the sides obtained by tracing the triangle in a clockwise direction.
Observing that $s_i$ is initial in $c$ if and only if the barring has $\upwide{i}$ and $\downwide{i+1}$, one easily verifies that property (i) of $c$-compatibility holds.

Conjugating~$c$ by an initial letter $s_i$ has the effect of swapping the positions of $i$ and $i+1$ in the cycle~$c$ and leaving the positions of the other entries unchanged.
Thus conjugating~$c$ by $s_i$ changes the circular diagram for~$c$ into the circular diagram for $scs$.
If~$W'$ is an irreducible rank-two parabolic, the triangle for~$W'$ in the circular diagram for~$c$ is located in the same place as
the triangle for $sW's$ in the circular diagram for $scs$.
Property (ii) is now easily checked, considering separately the cases $s\in W'$ and $s\not\in W'.$\

The orientation of a given~$W'$ can also be determined directly from the barring of $[2,n]$. 
As illustrated in Figures~\ref{cycleA}(b) and~\ref{cycleA}(c), if $j$ is upper-barred then~$W'$ is oriented 
\begin{equation}
\label{upj}
\big(\,\up{j\,}\,\,k\big)\toname{c}\big(i\,\,\,\up{j\,}\,\big)
\end{equation}
and if $j$ is lower-barred then~$W'$ is oriented 
\begin{equation}
\label{downj}
\big(i\,\,\down{\,j}\,\big)\toname{c}\big(\,\down{j}\,\,\,k\big).
\end{equation}

\begin{figure}[ht]
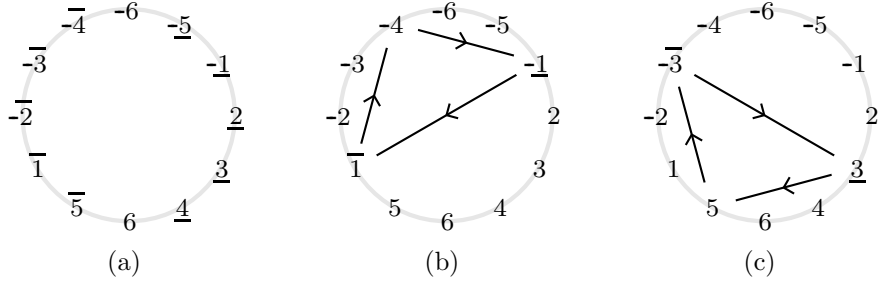

\centerline{
\begin{tabular}{ccccccc}
\scalebox{.8}{
\epsfbox{cycleBa.ps}
}
&&&
\scalebox{.8}{
\epsfbox{cycleBb.ps}
}
&&&
\scalebox{.8}{
\epsfbox{cycleBc.ps}
}
\\[4 pt]
(a)&&&(b)&&&(c)
\end{tabular}
}
\caption{$c$-Orientation of $B_6$}
\label{cycleB}
\end{figure}

We now consider the case $W=B_n$.
A Coxeter element~$c$ is realized as a signed permutation with a single cycle of length $2n$ of the form 
\[(-n\,\,d_1\,\,d_2\,\,\cdots\,\,d_l\,\,n\,\,u_k\,\,u_{k-1}\,\,\cdots\,\,u_1),\]
where $d_1,d_2,\ldots d_l$ are the lower-barred elements of $\pm[n-1]$ in increasing order and $u_1,u_2,\ldots,u_k$ are the upper-barred 
elements of $\pm[n-1]$ in increasing order.
As before, we place the elements of the cycle on a circle clockwise in the order indicated by the cycle, 
but we impose the extra condition that each $i$ be placed antipodally to $-i$.
Figure~\ref{cycleB}(a) illustrates for $c=s_5s_1s_0s_2s_3s_4$ in $B_6$ (cf. Figure~\ref{barB}).

Suppose two noncommuting reflections $t_1$ and $t_2$ are the canonical generators of the parabolic subgroup~$W'$ which contains them.
There are two possibilities.
One is that the two reflections are $(i\,\,j)$ and $(j\,\,k)$ for $i<j<k$ such that $|i|$, $|j|$ and $|k|$ are distinct.
In this case $|W'|=6$ and~$W'$ is represented as a centrally symmetric pair of inscribed triangles in the cycle representation of~$c$.
In each triangle the edges correspond to reflections in~$W'$ and we orient them cyclically clockwise.
In terms of the barring of $\pm[n-1]$, this orientation is given by (\ref{upj}) or (\ref{downj}). 

The other possibility is that the two reflections are $(-i\,\,i)$ and $(i\,\,j)$ for $-i<i<j$.
In this case $|W'|=8$ and~$W'$ permutes the symbols $\pm i,\pm j$, which are the vertices of an inscribed quadrilateral. 
The diagonal connecting $i$ to $-i$ splits the quadrilateral into two triangles.
The edges of either triangle correspond to the reflections $(i\,\,-i)$, $(i\,\,j)$ and $(-i\,\,j)$.
We orient~$W'$ by orienting these three reflections according to a clockwise walk around the triangle and inserting $(-j\,\,j)$ into the 
appropriate place in the cycle.
Specifically, there are two cases, depending on whether $i$ is upper-barred or lower-barred.
These are illustrated in Figures~\ref{cycleB}(b) and~\ref{cycleB}(c) and described below.
\begin{equation}
\label{upi}
\big(\,\up{i\,}\,\,\,j\big)\toname{c}\big(\,\down{-i\,}\,\,\,\upwide{i}\,\big),\mbox{ or}
\end{equation}
\begin{equation}
\label{downi}
\big(\,\up{-i\,}\,\,\,\downwide{i}\,\big)\toname{c}\big(\,\downwide{i}\,\,\,j\big).
\end{equation}
By an argument similar to the argument for type A, this family of orientations satisfies properties (i) and (ii) of $c$-compatibility.

\begin{figure}[ht]
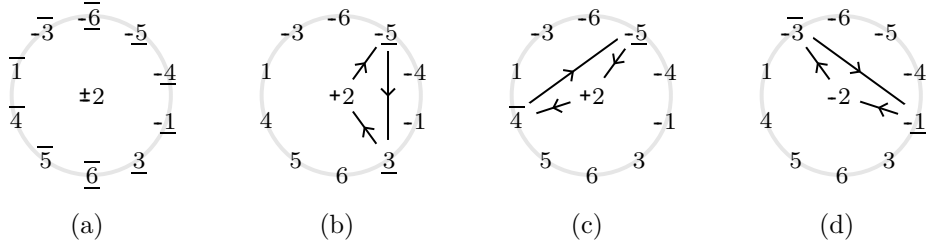

\centerline{
\begin{tabular}{cccc}
\scalebox{.75}{
\epsfbox{cycleDa.ps}
}
&
\scalebox{.75}{
\epsfbox{cycleDb.ps}
}
&
\scalebox{.75}{
\epsfbox{cycleDc.ps}
}
&
\scalebox{.75}{
\epsfbox{cycleDd.ps}
}
\\
(a)&(b)&(c)&(d)
\end{tabular}
}
\caption{$c$-Orientation of $D_6$}
\label{cycleD}
\end{figure}

We now consider the Coxeter group $D_n$.
A Coxeter element~$c$ has two cycles, one of length 2 and the other of length $2n-2$.
If~$c$ is symmetric then the short cycle is $(-1\,\,1)$ and if~$c$ is asymmetric then the short cycle is $(-2\,\,2)$.
In either case the long cycle is
\[(-n\,\,d_1\,\,d_2\,\,\cdots\,\,d_l\,\,n\,\,u_k\,\,u_{k-1}\,\,\cdots\,\,u_1),\]
where $d_1,d_2,\ldots d_l$ are the lower-barred elements of $\pm[2,n-1]$ or $\set{\pm 1}\cup\pm[3,n-1]$ in increasing order and 
$u_1,u_2,\ldots,u_k$ are the upper-barred elements in increasing order.
We place both elements of the short cycle at the center of a circle and place the elements of the large cycle on the circle clockwise in the 
order indicated by the cycle, with opposite elements placed antipodally, as illustrated in Figure~\ref{cycleD} for the asymmetric case $c=s_5s_4s_1s_2s_0s_3$ in $D_6$ (cf. Figure~\ref{barD}).
It is convenient to adopt the convention that $n$ is both upper-barred and lower-barred, and to adopt the same convention for $-n$.
The two elements at the center are called {\em central elements}, which are neither upper-barred nor lower-barred.

Two noncommuting reflections $t_1$ and $t_2$ are the canonical generators of the parabolic subgroup~$W'$ containing them if and only if they are $(i\,\,j)$ and $(j\,\,k)$ for $i<j<k$ with $|i|$, $|j|$ and $|k|$ distinct.
The subgroup~$W'$ is represented by the triangle with vertices $i$, $j$ and $k$, whose three sides are naturally identified with the 
reflections in~$W'$.
Orient~$W'$ so as to traverse these edges clockwise.
As in types A and B, one checks that this family of orientations satisfies properties (i) and (ii) of $c$-compatibility.

We now describe the orientations in terms of the barrings.
Suppose~$W'$ is represented by the triangle with vertices $i$, $j$ and $k$ for $i<j<k$ with $|i|$, $|j|$ and $|k|$ distinct.
In what follows, it is useful to keep in mind that~$W'$ is also represented by the triangle with vertices $-k$, $-j$ and $-i$.

If $j$ is not a central element then~$W'$ is oriented as in (\ref{upj}) or (\ref{downj}).
If $j$ is a central element then $j$ is neither upper-barred nor lower-barred.
The orientation of~$W'$ depends on the barring of $i$ and $k$.
If $i$ and $k$ have the same barring, then by possibly passing to the triangle with vertices $-k$, $-j$ and $-i$, we assume
that $i$ and $k$ are lower-barred and orient
\begin{equation}
\label{downidownk}
\big(j\,\,\,\down{\,k}\,\big)\toname{c}\big(\downwide{i}\,\,\,j\big).
\end{equation}
If $i$ and $k$ have opposite barrings then we assume $k<-i$ and orient:
\begin{equation}
\label{downiupk}
\big(\downwide{i}\,\,\,j\big)\toname{c}\big(j\,\,\,\up{k}\,\big),\mbox{ or}
\end{equation}
\begin{equation}
\label{upidownk}
\big(j\,\,\,\down{\,k}\,\big)\toname{c}\big(\,\up{i\,}\,\,\,j\big).
\end{equation}
To avoid confusion later, we again mention the convention that $n$ and $-n$ are upper-barred and lower-barred.
Figures~\ref{cycleD}(b), \ref{cycleD}(c) and~\ref{cycleD}(d) illustrate (\ref{downidownk}), (\ref{downiupk}) and (\ref{upidownk}) 
in the asymmetric case.

This construction of the $c$-orientation in types A, B and D completes the proof of Proposition~\ref{family}.
We conclude the section by discussing some additional properties of the $c$-orientation.

\begin{prop}
\label{restrict or}
For any $J\subseteq S$, the restriction of the $c$-orientation of~$W$ to the reflections in~$W_J$ coincides with the $c'$-orientation of~$W_J$,
where $c'$ is obtained from~$c$ by restriction.
\end{prop}
\begin{proof}
It is sufficient to treat the case where $J=\br{s}=S-\set{s}$ for some $s\in S$.
This is readily verified for the classical types via the description of the orientation in terms of planar diagrams, as follows:
Let~$W'$ be some irreducible component of~$W_{\br{s}}$ and let $c'$ be the Coxeter element of~$W'$ obtained from~$c$ by restriction.
One verifies that the planar diagram of $c'$ is obtained from the planar diagram for~$c$ by restriction, thus in particular preserving 
clockwise orientations of triangles.

The proposition is trivial for~$W$ of rank two and is checked for the exceptional types by computer, using the definition of the orientation
in terms of a special reduced word for~$w_0$.
%To do this in sort_maple, run:
%crys:=true: for R in [E6,E7,E8,F4] do check_res_lemma_all(R) od;
%crys:=false: for R in [H3,H4] do check_res_lemma_all(R) od;
\end{proof}

\begin{remark}\rm
\label{exten}
The orientation of the diagram of~$W$ induces a cyclic orientation of each standard irreducible rank-two parabolic subgroup of~$W.$\
In light of the case $|J|=2$ of Proposition~\ref{restrict or}, this orientation agrees with the $c$-orientation on these standard parabolic subgroups.
\end{remark}

\section{Aligned elements}
\label{align sec}
In this section we define $c$-aligned elements and show that they coincide with $c$-sortable elements.
As part of proving the coincidence, in types A, B and D we interpret the definition of $c$-aligned elements in terms of permutations.
These interpretations are also essential in proving crucial lemmas in Section~\ref{nc}.

An element $w\in W$ is {\em $c$-aligned} if for every irreducible rank-two parabolic subgroup~$W'$ with canonical generators 
$t_1$ and $t_2$ and $c$-orientation $t_1\toname{c} t_2$, if $t_2t_1t_2$ is an inversion of~$w$ then $t_1$ is an inversion of~$w$.
Equivalently, by Lemma~\ref{dyer lem}, the requirement is that if $I(w)\cap (W'-\set{t_2})$ is nonempty then $t_1$ is an inversion of~$w$.
The main result of this section is the following.
\begin{theorem}
\label{align}
An element~$w$ of~$W$ is $c$-sortable if and only if $w$ is $c$-aligned.
\end{theorem}

\begin{example}\rm
\label{B2align}
For $W=B_2$ with $c=s_0s_1$, the requirement that $w\in W$ be $c$-aligned is the following:
If $s_1s_0s_1$ is an inversion of~$w$ then $s_0$ is an inversion of~$w$.
The elements of~$W$ having $s_1s_0s_1$ as an inversion are $s_1s_0$, $s_1s_0s_1$, $s_0s_1s_0$ and $s_0s_1s_0s_1$.
Only the latter two have $s_0$ as an inversion, so the elements of~$W$ which fail to be $c$-aligned are $s_1s_0$ and $s_1s_0s_1$.
These are exactly the non-$c$-sortable elements (see Example~\ref{B2sort}).
\end{example}

\begin{remark}\rm
Billey and Braden~\cite{Bi-Br} developed a general theory of pattern avoidance for finite Coxeter groups and applied this theory to the 
study of Kazhdan-Lusztig polynomials.
The $c$-aligned elements constitute another instance of their theory.

A special case of the definition of the $c$-aligned elements is the definition of permutations with compressed inversion sets, used 
in \cite[Section~9]{Nonpure II} to study the Tamari lattice.
Compressed inversion sets were generalized in \cite[Section~5]{cambrian} and used to study Cambrian lattices of type A.
That generalization corresponds to the definition of $c$-aligned elements of type A.
In the crystallographic case, Ingalls and Thomas~\cite{IngTh} have arrived independently at a condition on inversion sets which is equivalent to $c$-alignment.
\end{remark}

Before proving Theorem~\ref{align}, we record two corollaries, the first of which holds because $I(w_0)=T$.
The second is immediate in light of Lemma~\ref{int}.

\begin{cor}
\label{w0 sort}
The longest element~$w_0$ of a finite Coxeter group~$W$ is $c$-sortable for any~$c$.
\end{cor}

\begin{cor}
\label{sort para}
Let $J\subseteq S$, let~$c$ be a Coxeter element of~$W$ and let $c'$ be the Coxeter element of~$W_J$ obtained by restriction.
If~$w$ is $c$-sortable then~$w_J$ is $c'$-sortable.
\end{cor}

The proof of Theorem~\ref{align} consists of showing that $c$-aligned elements share the recursive characterization of $c$-sortable elements described in Lemmas~\ref{sc} and~\ref{scs}.
Specifically, we prove the following lemmas.

\begin{lemma}
\label{align sc}
Let~$s$ be an initial letter of a Coxeter element~$c$ and let $w\in W$ have $l(sw)>l(w)$.
Then~$w$ is $c$-aligned if and only if it is an $sc$-aligned element of~$W_{\br{s}}$.
\end{lemma}

\begin{lemma}
\label{align scs}
Let~$s$ be an initial letter of a Coxeter element~$c$ and let $w\in W$ have $l(sw)<l(w)$.
Then~$w$ is $c$-aligned if and only if $sw$ is $scs$-aligned.
\end{lemma}

As a first step towards proving Lemmas~\ref{align sc} and~\ref{align scs}, we characterize the aligned elements in the classical types in terms of permutations.
When Theorem~\ref{align} is proved, these characterizations will apply to sortable elements (see Theorem~\ref{ABDsort}).

A permutation $\pi\in S_{n+1}$ satisfies condition (A) if both of the following conditions hold:
\begin{enumerate}
\item[(A1) ]$\pi$ contains no subsequence $\upwide{j}k\,i$ with $i<j<k$, and
\item[(A2) ]$\pi$ contains no subsequence $k\,i\downwide{j}$ with $i<j<k$.
\end{enumerate}
A ``subsequence'' of~$\pi$ means a subsequence of the one-line notation for~$\pi$.

\begin{lemma}
\label{Aalign}
A permutation $\pi\in S_{n+1}$ is $c$-aligned if and only if it satisfies condition (A) with respect to the barring corresponding to~$c$.
\end{lemma}

\begin{proof}
In light of (\ref{upj}) and (\ref{downj}),~$\pi$ is $c$-aligned if and only if for each $i<j<k$ the following conditions hold:
If $(i\,\,k)\in I(\pi)$ and $j$ is upper-barred then $\big(\,\up{j\,}\,\,\,k\big)\in I(\pi)$ and if $(i\,\,k)\in I(\pi)$ and $j$ is lower-barred then 
$\big(i\,\,\,\down{\,j}\,\big)\in I(\pi)$.
These conditions are equivalent to (A1) and (A2).
\end{proof}

A signed permutation $\pi\in B_n$ satisfies condition (B) if the long one-line notation for~$\pi$ has no subsequence 
$\upwide{j}k\,i$ with $i<j<k$.
Equivalently,~$\pi$ satisfies (B) if it has no subsequence $k\,i\downwide{j}$ with $i<j<k$.

\begin{lemma}
\label{Balign}
A signed permutation $\pi\in B_n$ is $c$-aligned if and only if it satisfies condition (B) with respect to the barring corresponding to~$c$.
\end{lemma}

\begin{proof}
The alignment condition on rank-two parabolic subgroups with $|W'|=6$ becomes the requirement that~$\pi$ has no subsequence 
$\upwide{j}k\,i$ and equivalently has no subsequence $k\,i\downwide{j}$ for $i<j<k$ with $i$, $j$ and $k$ having distinct absolute values. 
In light of (\ref{upi}) and (\ref{downi}), the alignment condition on rank-two parabolic subgroups with $|W'|=8$ is equivalent to the 
requirement that~$\pi$ have no subsequence of any of the following forms for $-i<i<j$:
\[\upwide{i}-\!j\,j\,\down{-i\,},\quad  \upwide{i}\,j-\!j\,\down{-i\,},\quad  \upwide{-i}\,j-\!j\,\downwide{i}\quad \mbox{or}\quad
j\,\up{-i\,}\,\downwide{i}-\!j.\]
This is equivalent to the requirement that~$\pi$ have no subsequence $\upwide{j}k\,i$ with $i<j<k$ such that $i$, $j$ and $k$ do not have three distinct absolute values.
\end{proof}

\begin{example}\rm
We continue Example~\ref{B2align}.
The barring corresponding to $c$ has $\upwide{-1}$.
Thus condition (B) requires that the long one-line notation contain neither the subsequence $(-1)\,2\,(-2)$ nor the subsequence $(-1)\,1\,(-2)$.
The element $s_1s_0=(-1)\,2\,(-2)\,1$ contains the former subsequence and $s_1s_0s_1=2\,(-1)\,1\,(-2)$ contains the latter.
The other six elements of $B_2$ contain neither subsequence.
\end{example}

An even-signed permutations~$\pi$ satisfies condition (D) if and only if all four of the following conditions hold: 
\begin{enumerate}
\item[(D1) ]$\pi$ has no subsequence $\upwide{j}\,k\,i$ with $i<j<k$ such that $|i|$, $|j|$ and $|k|$ are distinct,
\item[(D2) ]$\pi$ has no subsequence $j\,\down{k}\,\down{\,i}$ with $i<j<k$ such that $|i|$, $|j|$ and $|k|$ are distinct and $j$ is central.
\item[(D3) ]$\pi$ has no subsequence $j\,\up{k}\,\down{\,i}$ with $-k<i<j<k$ such that $|i|$, $|j|$ and $|k|$ are distinct and $j$ is central.
\item[(D4) ]$\pi$ has no subsequence $j\,\down{k}\,\up{\,i\,}$ with $i<j<k<-i$ such that $|i|$, $|j|$ and $|k|$ are distinct and $j$ is central.
\end{enumerate}
In interpreting condition (D), recall that $n$ and $-n$ are both upper-barred and lower-barred and that central elements are neither 
upper-barred nor lower-barred.
Note also that for each of these conditions there is an equivalent formulation obtained by applying the symmetry of the long one-line notation for~$\pi$.

\begin{lemma}
\label{Dalign}
An even-signed permutation $\pi\in D_n$ is $c$-aligned if and only if it satisfies condition (D) with respect to the barring 
corresponding to~$c$.
\end{lemma}

\begin{proof}
Condition (D1) arises from rank-two parabolic subgroups~$W'$ corresponding to the case where $j$ is not central.
Conditions (D2) and (D4) come from the orientations described in (\ref{downidownk}) and (\ref{upidownk}).
Condition (D3) is obtained from the orientation described in (\ref{downiupk}) by the symmetry of the long one-line notation.
\end{proof}

We now prove Lemmas~\ref{align sc} and~\ref{align scs}.

\begin{proof}[Proof of Lemma~\ref{align sc}]
Let~$w$ be an $sc$-aligned element of~$W_{\br{s}}$.
Then by Lemma~\ref{int} and Proposition~\ref{restrict or},~$w$ is $c$-aligned as well.
Conversely, suppose~$w$ is $c$-sortable.
If $w\in W_{\br{s}}$ then by Lemma~\ref{int} and Proposition~\ref{restrict or},~$w$ is $sc$-aligned.
Thus it remains to prove that $w\in W_{\br{s}}$.

In the classical types, we represent~$w$ by a permutation~$\pi$ and suppose $\pi\not\in W_{\br{s}}$.
We will reach a contradiction to Lemma~\ref{Aalign}, \ref{Balign}, or~\ref{Dalign}.

If $W=S_{n+1}$ and $s=s_j$, then $\upwide{j}$ precedes $\downwide{j+1}$.
Because $\pi\not\in W_{\br{s}}$, some entry $k\ge j+1$ precedes some entry $i\le j$ in the one-line notation for~$\pi$.
If $k$ follows $j$ then $\upwide{j}k\,i$ is a subsequence of~$\pi$ violating (A1) and if $k$ precedes $j$ then $k\,j\downwide{j+1}$ violates (A2).

Now suppose $W=B_n$.
If $s=s_j$ for some $j\in[n-1]$ then let $k=j+1$.
If $s=s_0$ then let $j=-1$ and $k=1$.
Then $\upwide{j}$ precedes $\downwide{k}$ and some entry $l\ge k$ precedes some entry $i\le j$ in the long one-line notation for~$\pi$.
If $l$ follows $j$ then $\upwide{j}l\,i$ is a violation of (B) in~$\pi$ and if $l$ precedes $j$ then $l\,j\downwide{k}$ is a violation of (B).

Next we show that the lemma for $D_n$ reduces to the case $n\le 6$.
Let $W$ be of type $D_n$ and let $(i\,\,\,j)$ be a transposition associated to $s$.
Since $\pi\not\in W_{\br{s}}$, there is some transposition $(k\,\,\,l)$ corresponding to an inversion of~$\pi$ not in $W_{\br{s}}$.
As explained in Remark~\ref{Dn res}, let~$\pi'$ be the even-signed permutation obtained by restricting~$\pi$ to the set $\set{\pm 1,\pm 2,\pm i,\pm j,\pm k,\pm l}$.
Let $W'$, $n'$ and $c'$ be as in Remark~\ref{Dn res}, so that in particular $n'\le 6$.
The simple generator $s$ is also a simple generator of $W'$, initial in $c'$  with $l(s\pi')<l(\pi')$.
The other simple generators of $W'$ are all in particular elements of $W_{\br{s}}$, so the fact that $(k\,\,\,l)\not\in W_{\br{s}}$ implies that $(k\,\,\,l)\not\in W'_{\br{s}}$.
The transposition $(k\,\,\,l)$ is an inversion of~$\pi'$, so $\pi'\not\in W'_{\br{s}}$.
By the case $W=D_n$ for $n\le 6$,~$\pi'$ violates (D) with respect to~$c'$, and the three entries of~$\pi$ constituting this violation of (D) in~$\pi'$ also constitute a violation of (D) in~$\pi$.

In rank two, in the exceptional cases and for $D_n$ for $n\le 6$, the $c$-orientation was constructed in terms of a $c$-sorting word for~$w_0$.
We further verify the following claim.
Let~$s$ be initial in~$c$ and let $t\in T$ have $t\not\in W_{\br{s}}$ and $t\neq s$.
We claim that there exists $t'\not\in W_{\br{s}}$ with the following two properties:
$t'$ precedes $t$ in the reflection order arising from the $c$-sorting word for~$w_0$ and
the rank-two parabolic subgroup containing $t$ and $t'$ has canonical generators $t'$ and $t''$ for some $t''\neq t$.
This claim is trivial in rank two and is checked computationally for $D_n$ with $n\le 6$ and the exceptional types.
%To do this in sort_maple, run:
%crys:=true: for R in [E6,E7,E8,F4] do check_claim_all(R) od;
%crys:=false: for R in [H3,H4] do check_claim_all(R) od;

We use this claim to verify the lemma in the remaining cases.
Suppose $t$ is a reflection not in~$W_{\br{s}}$.  
Choose a reduced word for~$c$ in which~$s$ is the first letter, and use this reduced word for~$c$ to construct a $c$-sorting word for~$w_0$.
We show, by induction on the reflection order arising from this $c$-sorting word, that $t$ is not an inversion of~$w$.
If $t$ is first in the reflection order then $t=s$ and by hypothesis~$s$ is not an inversion of~$w$.
If $t$ comes later then let $t'$ be as in the claim.
By induction, $t'$ is not an inversion of~$w$, and therefore since~$w$ is $c$-aligned, $t$ is not an inversion of~$w$.
We have shown that every inversion of~$w$ is in~$W_{\br{s}}$ and thus that~$w$ itself is in~$W_{\br{s}}$.
\end{proof}

\begin{proof}[Proof of Lemma~\ref{align scs}]
Let~$W'$ be an irreducible parabolic subgroup of rank two with canonical generators $t_1$ and $t_2$.
We claim that~$W'$ is an exception to the $c$-alignment of~$w$ if and only if $sW's$ is an exception to the $scs$-alignment of $sw$.
If $t_1,t_2\neq s$ then the claim is immediate by property (ii) of $c$-compatibility and the fact that $I(sw)=sI(w)s-\set{s}$.
Otherwise, take $t_1=s$, so that $s\toname{c}t_2$ and $t_2\toname{scs}s$.
In this case~$W'$ is not an exception to the $c$-alignment of~$w$ because $s\in I(w)$.
If $sW's=W'$ is an exception to the $scs$-alignment of $sw$ then $st_2s\in I(sw)$ but $t_2\not\in I(sw)$.
But then by Lemma~\ref{dyer lem}, $s\in I(sw)$, contradicting the hypotheses on~$w$.
This contradiction shows that $sW's$ is not an exception to the $scs$-alignment of $sw$, completing the proof of the claim.
Letting~$W'$ vary over all irreducible rank-two parabolics proves the lemma.
\end{proof}

In light of Lemmas~\ref{align sc} and~\ref{align scs}, the proof of Theorem~\ref{align} is simple.
Given $w\in W$ and~$s$ initial in $c$, if $l(sw)>l(w)$ then one uses Lemmas~\ref{sc} and~\ref{align sc} to argue  by induction on the rank of~$W.$\
If $l(sw)<l(w)$ then one uses Lemmas~\ref{scs} and~\ref{align scs} to argue  by induction on length.
The base of the induction is the fact that $1$ is both $c$-sortable and $c$-aligned for any $c$.

Theorem~\ref{align} and Lemmas~\ref{Aalign}, \ref{Balign} and~\ref{Dalign} now imply the following:

\begin{theorem}
\label{ABDsort}
For $X\in\set{A,B,D}$, let $W$ be a finite Coxeter group of type~$X$.
An element $w\in W$ is $c$-sortable if and only if the associated permutation satisfies condition (X) with respect to the barring corresponding to $c$.
\end{theorem}

\section{Noncrossing partitions}
\label{nc prelim}
In this section we present background on noncrossing partitions based on \mbox{\cite{Bessis,BWorth}}.

Any element $w\in W$ can be written as a word in the alphabet $T$.
To avoid confusion we always refer to a word in the alphabet $T$ as a \emph{$T$-word}.
Any other use of the term ``word'' is assumed to refer to a word in the alphabet~$S$.
A \emph{reduced} $T$-word for~$w$ is a $T$-word for~$w$ which has minimal length among all $T$-words for~$w$.
The \emph{absolute length} of an element~$w$ of~$W$ is the length of a reduced $T$-word for~$w$.
This is not the usual length $l(w)$ of~$w$, which is the length of a reduced word for~$w$ in the alphabet~$S$.

The notion of reduced $T$-words leads to a prefix partial order on~$W\!$, analogous to the weak order.
Say $x\le_T y$ if~$x$ possesses a reduced $T$-word which is a prefix of some reduced $T$-word for~$y$.
Equivalently, $x\le_Ty$ if every reduced $T$-word for~$x$ is a prefix of some reduced $T$-word for~$y$.
Since the alphabet $T$ is closed under conjugation by arbitrary elements of~$W\!$, the partial order $\le_T$ is invariant under conjugation.
The partial order $\le_T$ can also be defined as a subword order:
$x\le_Ty$ if and only if there is a reduced $T$-word for $y$ having as a subword some reduced $T$-word for~$x$.
In particular, $x\le_Ty$ if and only if $x^{-1}y\le y$.

The {\em noncrossing partition lattice} in~$W$ (with respect to the Coxeter element~$c$) is the interval $[1,c]_T$, and the 
elements of this interval are called {\em noncrossing partitions}.
The poset $[1,c]_T$ is graded and the rank of a noncrossing partition is its absolute length.

The following theorem is \cite[Theorems 1 and 2]{BWorth}, rephrased in light of \cite[Proposition~1]{BWorth}.
For $x\in W$, let $U_x$ be the subspace fixed by~$x$, where~$x$ is identified with an orthogonal transformation via the reflection representation of~$W$ discussed near the beginning of Section~\ref{cox sec}.
\begin{theorem}
\label{subspaces}
If $x\le_T c$ and $y\le_T c$ then $x\le_T y$ if and only if $U_x\supseteq U_y$.
In particular, $U_x=U_y$ if and only if $x=y$.
\end{theorem}

If $t_1t_2\cdots t_l$ is any reduced $T$-word for $x\le_Tc$ then $U_x$ is equal to the intersection $\cap H_i$ for $i\in[l]$, where $H_i=U_{t_i}$ is the reflecting hyperplane associated to $t_i$.
The usual bijection from intersections of reflecting hyperplanes to parabolic subgroups takes $U_x$ to the parabolic subgroup $\br{t_1,t_2,\ldots,t_l}$.
Thus in particular, for each noncrossing partition there is a corresponding parabolic subgroup.
The noncrossing partition lattice is reverse containment on the parabolic subgroups arising in this manner.
(Except in ranks 1 or 2, this is a proper subset of the set of parabolic subgroups of~$W.$)
In particular, any reflection $t$ in $\br{t_1,t_2,\ldots,t_l}$ has $t\le_Tx$.

The following lemma summarizes the behavior of noncrossing partitions with respect to standard parabolic subgroups.
\begin{lemma}
\label{abs para}
If $x\le_Tc$ and $s\in S$ then the following are equivalent:
\begin{enumerate}
\item[(i) ]$x\in W_{\br{s}}$.
\item[(ii) ]Every reflection $t$ in any reduced $T$-word for~$x$ has $t\in W_{\br{s}}$.
\end{enumerate}
Under the additional condition that~$s$ is initial in~$c$, the following are equivalent to (i) and (ii):
\begin{enumerate}
\item[(iii) ]$x\le_Tsc$.
\item[(iv) ]$x\le_Tsx\le_Tc$.
\end{enumerate}
\end{lemma}
\begin{proof}
Recall that $w\in W_{\br{s}}$ if and only if~$w$ fixes $\rho(s)$, or in other words if and only if $\rho(s)\subseteq U_w$.
(Here, $\rho(s)$ is the ray defined just before Lemma~\ref{cover para}).
Since $U_x$ is the intersection $\cap U_{t_i}$ for any reduced $T$-word $t_1\cdots t_l$ for~$x$, conditions (i) and (ii) are equivalent.

Now suppose~$s$ is initial in $c$ so that in particular $sc\le_Tc$.
The subspace $U_{sc}$ is the line containing $\rho(s)$, so Theorem~\ref{subspaces} says that conditions (i) and (iii) are equivalent.

We now show that (iv) is equivalent to the previous three conditions.
If (iv) holds then there is a reduced $T$-word $st_1t_2\cdots t_{n-1}$ for~$c$ such that $t_1t_2\cdots t_k$ is a reduced $T$-word for~$x$ for some $k\le n-1$.
Therefore~$x$ is a prefix of the reduced $T$-word $t_1t_2\cdots t_{n-1}$ for $sc$ and thus $x\le_T sc$, so (i) holds.

Conversely, suppose (ii) and (iii).
By (ii), no reduced $T$-word for~$x$ begins with~$s$.
Thus prepending~$s$ to any reduced $T$-word for~$x$ gives~a reduced $T$-word and the subword characterization of $\le_T$ says that $x\le_Tsx$.
By (iii), there is a reduced $T$-word for~$x$ which is a prefix of some reduced $T$-word for $sc$.
Prepending~$s$ gives~a reduced $T$-word for $sx$ which is a prefix of some reduced $T$-word for~$c$, so that $sx\le_Tc$.
\end{proof}

\section{Sortable elements and noncrossing partitions}
\label{nc sec}
In this section, we define a bijection between sortable elements and noncrossing partitions.
Let~$w$ be a $c$-sortable element and let $a=a_1a_2\cdots a_k$ be a $c$-sorting word for~$w$.
Totally order the inversions of~$w$ such that the $i$th reflection in the order is $a_1a_2\cdots a_{i-1}a_ia_{i-1}\cdots a_2a_1$.
Equivalently, $t$ is the $i$th reflection in the order if and only if $tw=a_1a_2\cdots \hat{a}_i\cdots a_k$, where $\hat{a}_i$ indicates that $a_i$ is deleted from the word. 
Write the set of cover reflections of~$w$ as a subsequence $t_1,t_2,\ldots,t_l$ of this order on inversions.
Let $\nc_c$ be the map which sends~$w$ to the product $t_1t_2\cdots t_l$.
Recall that the construction of a $c$-sorting word begins with an arbitrary choice of a reduced word for~$c$.
However, since any two $c$-sorting words for~$w$ are related by commutation of simple generators,
$\nc_c(w)$ does not depend of the choice of reduced word for~$c$.

\begin{theorem}
\label{nc}
For any Coxeter element~$c$, the map $w\mapsto\nc_c(w)$ is a bijection from the set of $c$-sortable elements to 
the set of noncrossing partitions with respect to~$c$.
Furthermore $\nc_c$ maps $c$-sortable elements with $k$ descents to noncrossing partitions of rank $k$.
\end{theorem}
Recall that the descents of~$w$ are the simple generators $s\in S$ such that $l(ws)<l(w)$.
Recall also that these are in bijection with the cover reflections of~$w$.

\begin{example}\rm
\label{B2nc}
We again consider the case $W=B_2$ and $c=s_0s_1$.
As a special case of the combinatorial realization of noncrossing partitions of type B given in~\cite{Reiner}, the noncrossing partitions in $B_2$ with respect to~$c$ are the centrally symmetric noncrossing partitions of the cycle shown below.
\[\scalebox{.8}{\epsfbox{B2cycle.ps}}\]
Figure~\ref{B2ncfig} illustrates the map $\nc_c$ for this choice of~$W$ and~$c$.
\end{example}

\begin{figure}[ht]
\[\begin{array}{|c||c|c|c|c|c|c|}\hline
&&&&&&\\[-7 pt]
w&	1&	\hat{s}_0&	s_0\hat{s}_1&	s_0s_1|\hat{s}_0&	\hat{s}_0s_1|s_0\hat{s}_1&	\hat{s}_1\\[3 pt]\hline
&&&&&&\\[-7 pt]
nc_c(w)&	1&	s_0&		s_0s_1s_0&	s_1s_0s_1&	s_0\cdot s_1&		s_1\\[3 pt]\hline
&&&&&&\\[-7 pt]
&	\epsfbox{B2.1.ps}&	\epsfbox{B2.a.ps}&\epsfbox{B2.ab.ps}&\epsfbox{B2.aba.ps}&\epsfbox{B2.ababi.ps}&\epsfbox{B2.b.ps}\\
&&&&&\downarrow&\\
&&&&&\epsfbox{B2.ababii.ps}&\\
\hline
\end{array}\]
\caption{The map $\nc_c$}
\label{B2ncfig}
\end{figure}

\begin{example}\rm
\label{Annc}
Covering reflections of a permutation $\pi\in S_{n+1}$ are the transpositions corresponding to descents (pairs $(\pi_i,\pi_{i+1})$ with $\pi_i>\pi_{i+1}$).
The map $\nc_c$ sends~$\pi$ to the product of these transpositions.
The relations $\pi_i\equiv\pi_{i+1}$ for descents $(\pi_i,\pi_{i+1})$ generate an equivalence relation on $[n+1]$ which can be interpreted as a noncrossing partition (in the classical sense)
of the cycle~$c$.
For $c=(n\,\,\,\,n\!+\!1)\cdots(2\,\,3)(1\,\,2)$ as in Example~\ref{Ansort}, this map between $231$-avoiding permutations 
and classical (i.e.\ type A) noncrossing partitions is presumably known.
\end{example}

\begin{remark}\rm
\label{inf nc}
The definition of $\nc_c$ is valid for infinite $W$.
However, the proofs presented in this section address the finite case only.
In particular, for infinite $W$ it is not even known whether $\nc_c$ maps $c$-sortable elements into the interval $[1,c]_T$.
\end {remark}

The basic tool for proving Theorem~\ref{nc} is induction on rank and length, using Lemmas~\ref{sc} and~\ref{scs}.
A more complicated induction is used to prove the existence of the inverse map.
The inductive argument requires the following three lemmas.

\begin{lemma}
\label{nc s}
Let~$s$ be initial in~$c$ and let~$w$ be $c$-sortable with $l(sw)<l(w)$.
If~$s$ is a cover reflection of~$w$ then $\nc_{scs}(sw)=\nc_c(w)\cdot s$.
If~$s$ is not a cover reflection of~$w$ then $\nc_{scs}(sw)=s\cdot\nc_c(w)\cdot s$.
\end{lemma}
\begin{proof}
Use a reduced word for $c$ whose first letter is $s$ to write a $c$-sorting word for $w$.
Delete the first letter from this word to obtain an $scs$-sorting word for $sw$.
For these choices of sorting words, compare the defining factorization of $\nc_c(w)$ to the defining factorization of $\nc_{scs}(sw)$.
The lemma is immediate.
\end{proof}

\begin{lemma}
\label{nc lemma}
Let $s\in S$ be either initial or final in~$c$ and let~$w$ be $c$-sortable.
If~$s$ is a cover reflection of~$w$ then every other cover reflection of~$w$ is in~$W_{\br{s}}$.
\end{lemma}
\begin{proof}
For~$W$ of rank two, the lemma is trivial. 
In the exceptional cases and in the case $D_n$ with $n\le 6$, it is checked by computer.
%To do this in sort_maple, run:
%crys:=true: for R in [A3,B3,B4,D4,E6,E7,E8,F4] do check_nc_lemma_all(R) od;
%crys:=false: for R in [H3,H4] do check_nc_lemma_all(R) od;
For classical~$W,$\ let~$\pi$ be the permutation corresponding to~$w$.

First consider the case $W=S_{n+1}$.
Suppose~$\pi$ satisfies (A) with respect to $c$ and let $s_j$ be initial or final in $c$ such that $s_j$ is a cover reflection of~$\pi$.
Then $j+1$ immediately precedes $j$, and one of the two is upper-barred while the other is lower-barred.
Let $(i\,\,\,k)$ be any other cover reflection (with $i<k$).
Since $j+1$ immediately precedes $j$ and $k$ immediately precedes $i$, either $k$ follows (or coincides with) $j$ or $i$ precedes 
(or coincides with) $j+1$.
In either case, $i\neq j$ and $k\neq j+1$.
In the former case, if $k>j+1>j>i$ then~$\pi$ contains either $\upwide{j}k\,i$ or $\upwide{j+1}ki$, violating (A1).
In the latter case, if $k>j+1>j>i$ then~$\pi$ contains either $k\,i\downwide{j}$ or $k\,i\downwide{j+1}$, violating (A2).
Thus either $i\ge j+1$ or $k\le j$, so that $(i\,\,j)\in W_{\br{s_j}}$.
The proof for the case $W=B_n$ is nearly identical.

Finally, consider $W$ of type $D_n$ and let $(i\,\,\,j)$ be a transposition associated to $s$.
Let $(k\,\,\,l)$ be another cover reflection of~$\pi$.
We want to show that $(k\,\,\,l)\in W_{\br{s}}$.
As explained in Remark~\ref{Dn res}, let~$\pi'$ be the signed permutation obtained by restricting~$\pi$ to the set $\set{\pm 1,\pm 2,\pm i,\pm j,\pm k,\pm l}$.
For $W'$, $c'$ and $n'$ as in Remark~\ref{Dn res}, $s$ is a simple generator of $W'$ which is either initial or final in $c'$.
Now $(k\,\,\,l)$ is an inversion of~$\pi'$, and because~$\pi$ satisfies (D) with respect to $c$,~$\pi'$ satisfies (D) with respect to $c'$.
Thus by the case $n\le 6$, $(k\,\,\,l)\in W'_{\br{s}}$.
The simple generators of $W'$, besides $s$, are all elements of $W_{\br{s}}$, so $(k\,\,\,l)\in W_{\br{s}}$.
\end{proof}

\begin{lemma}
\label{nc lemma 2}
Let~$s$ be final in~$c$ and let $w\in W$ be $c$-sortable.
Then $l(sw)<l(w)$ if and only if~$s$ is a cover reflection of~$w$.
\end{lemma}
\begin{proof}
The ``if'' statement is by the definition of cover reflection, so we suppose $l(sw)<l(w)$ and prove that $s$ is a cover reflection.
The rank-two case is trivial and the exceptional cases are checked by computer.
%To do this in sort_maple, run:
%crys:=true: for R in [E6,E7,E8,F4] do check_nc_lemma2_all(R) od;
%crys:=false: for R in [H3,H4] do check_nc_lemma2_all(R) od;
In each classical case, let~$\pi$ be the permutation corresponding to~$w$.

First suppose $W=S_{n+1}$ and suppose $s=s_i$.
Then since~$s$ is final and $l(s\pi)<l(\pi)$, $\upwide{i+1}$ precedes $\downwide{i}$ in the one-line notation for~$\pi$.
If any other entry $k$ occurs between them then the sequence $\upwide{i+1}k\downwide{i}$ violates either (A1) or (A2).
Thus they are adjacent, so that $s_i$ is a cover reflection of~$\pi$.
The case $W=B_n$ is proved similarly.

For $W=D_n$, suppose $s=s_i$ for $i\ge 2$.
Then $\upwide{i+1}$ precedes $\downwide{i}$ in the one-line notation for~$\pi$.
Arguing as in type A we see that if any entry besides $-i$ or $-i-1$ occurs between $i$ and $i+1$ then~$\pi$ violates (D1), perhaps applying the symmetry of the long one-line notation.
Thus $s$ is a cover reflection for~$\pi$.
If $s=s_0$ or $s=s_1$ then by symmetry we assume that $s=s_1$.
Thus $2$ precedes $1$ in~$\pi$ and since~$s$ is final, either $1$ is lower-barred and $2$ is central or $2$ is upper-barred and $1$ is central.
We want to show that no element, except perhaps $-1$ or $-2$ appears between $2$ and $1$.
Suppose $k>2$ appears between $2$ and $1$.
If the barring has $\downwide{1}$ then either $2\downwide{k}\,\downwide{1}$ violates (D2) or $2\upwide{k}\downwide{1}$ violates (D3).
If the barring has $\upwide{2}$ then $\upwide{2}k\,1$ violates (D1).
Suppose $k<-2$ appears between $2$ and $1$, so that by symmetry $-1\,-k\,-2$ is a subsequence of~$\pi$.
If the barring has $\downwide{1}$ then $\upwide{-1}-k\,-2$ violates (D1).
If the barring has $\upwide{2}$ then either $-1\downwide{-k}\,\downwide{-2}$ violates (D2) or $-1\upwide{-k}\downwide{-2}$ violates (D3).
\end{proof}

We close the section with a proof of the main theorem followed by a remark.

\begin{proof}[Proof of Theorem~\ref{nc}]
First we prove that for a $c$-sortable element~$w$, the $T$-word $t_1t_2\cdots t_l$ arising in the definition of $\nc_c(w)$ is a reduced 
$T$-word for some $x\le_Tc$.
This shows that $\nc_c$ indeed maps $c$-sortable elements to noncrossing partitions with respect to~$c$ and also proves the second
statement in the theorem.

Let~$s$ be the first letter of the reduced word for $c$ used in defining $\nc_c$.
If $l(sw)>l(w)$ then by Lemma~\ref{sc},~$w$ is an $sc$-sortable element of~$W_{\br{s}}$.
By induction on rank, the $T$-word arising in the definition of $\nc_{sc}(w)$ is a reduced $T$-word for some element $x\le_Tsc$.
But this $T$-word is identical to the $T$-word arising in the definition of $\nc_c(w)$, and since $sc\le_Tc$, we are finished in the case where
$l(sw)>l(w)$.

If $l(sw)<l(w)$, we consider two cases.
If~$s$ is not a cover reflection of~$w$ then the $T$-word arising from the definition of $\nc_{scs}(sw)$ is $t'_1\cdots t'_l$, where $t'_i=st_is$ for
each $i\in[l]$.
By Lemma~\ref{scs} and  induction on $l(w)$, this is a reduced $T$-word for an element $y\le_Tscs$.
Conjugating by~$s$, we see that $t_1\cdots t_l$ is a reduced $T$-word for an element below~$c$.

If~$s$ is a cover reflection of~$w$ then in particular $s=t_1$.
In this case, the $T$-word arising from the definition of $\nc_{scs}(sw)$ is $t'_2\cdots t'_l$, with $t'_i$ as in the previous paragraph.
By induction on the length of~$w$, this is a reduced $T$-word for an element $y\le_Tscs$.
Thus $t_2\cdots t_l$ is a reduced $T$-word for an element $x=sys\le_Tc$.
Lemmas~\ref{cover para} and~\ref{nc lemma} imply that $t_2\cdots t_l\in W_{\br{s}}$, so by Lemma~\ref{abs para}, $st_2\cdots t_l$ is a reduced $T$-word for $sx$ and $sx\le_Tc$.
We have shown that $t_1\cdots t_l$ is a reduced $T$-word for a noncrossing partition.

We now show that for any~$x$ with $x\le_Tc$, there exists a unique $c$-sortable element mapped to~$x$ by $\nc_c$.
The argument is by induction on the rank of~$W,$\ on the absolute length of~$x$ and on $\nu_c(W')$, where~$W'$ is the parabolic subgroup associated to~$x$ and $\nu_c$ is as defined in the paragraph following Lemma~\ref{nu lemma}.

If $x\in W_{\br{s}}$ for some initial letter~$s$ of~$c$ then by Lemma~\ref{abs para}, $x\le_Tsc$.
By induction on the rank of~$W,$\ there is a unique $sc$-sortable element~$w$ of~$W_{\br{s}}$ with $\nc_{sc}(w)=x$.
The element~$w$ is $c$-sortable with $\nc_c(w)=x$ as well.
Any other $c$-sortable element~$w'$ of~$W$ with $\nc_c(w')=x$ has $w'\in W_{\br{s}}$ by Lemmas~\ref{cover para} and~\ref{abs para}.
Thus~$w$ is the unique $c$-sortable element in all of~$W$ with $\nc_c(w)=x$.

Now suppose there is no initial letter~$s$ of~$c$ such that $x\in W_{\br{s}}$.
In particular $x\neq 1$ and the corresponding parabolic subgroup~$W'$ is nontrivial.
Let $S'$ be the set of canonical generators of~$W'$.

If $\nu_c(W')=0$ then there is some initial letter~$s$ of~$c$ with $s\in W'$.
Then $s\le_Tx$ so that $sx<_Tx\le_Tc$.
Conjugating by~$s$ we have $xs\le_Tscs$.
By induction on absolute length there is a unique $scs$-sortable element $w\in W$ with 
$\nc_{scs}(w)=xs$.
If $l(sw)<l(w)$ then by Lemma~\ref{nc lemma 2},~$s$ is a cover reflection of~$w$.
However, this would mean in particular that there is a reduced $T$-word for $xs$ containing the letter~$s$, contradicting the fact that $xs<_Tx$.
Thus $l(sw)>l(w)$ and by Lemma~\ref{scs}, $sw$ is $c$-sortable.
We now claim that~$s$ is a cover reflection of $sw$.
If not then Lemma~\ref{nc s} says $\nc_c(sw)=s\cdot\nc_{scs}(w)\cdot s=sx\le_Tx$.
By Lemma~\ref{abs para}, $sx$ is in~$W_{\br{s}}$, so by Lemmas~\ref{cover para} and~\ref{abs para}, $sw\in W_{\br{s}}$, contradicting the fact that $l(s\cdot sw)<l(sw)$.
Therefore~$s$ is a cover reflection of $sw$ and by Lemma~\ref{nc s}, $\nc_c(sw)=\nc_{scs}(w)\cdot s=x$.

If~$w'$ is any other $c$-sortable element of~$W$ with $\nc_c(w')=x$ then since $x\not\in W_{\br{s}}$, $w'\not\in W_{\br{s}}$ and therefore $l(sw')<l(w')$.
Since $s\le_T x$,~$s$ is in $W'$.
By Lemma~\ref{s in},~$s$ is in $S'$ and therefore it is a cover reflection of~$w'$ by Lemma~\ref{canon lem}.
Thus by Lemma~\ref{nc s}, $\nc_{scs}(sw')=\nc_c(w')\cdot s=xs$.
By the uniqueness of $w$, we have $sw'=w$, so $sw$ is the unique $c$-sortable element mapping to~$x$ under $\nc_c$.

If $\nu_c(W')>0$ then let $s=s_0$ in some sequence $s_0,\ldots,s_{\nu_c(W')}$ satisfying the conclusion of Lemma~\ref{nu lemma}.
In particular~$s$ is initial in~$c$ and by Lemma~\ref{canon conj}, $sS's$ is the set of canonical generators for $sW's$, the parabolic subgroup corresponding to $y=sxs\le_Tscs$.
Since $\nu_{scs}(sW's)=\nu_c(W')-1$, by induction there is a unique $scs$-sortable element~$w$ with \mbox{$\nc_{scs}(w)=y$}.
If $l(sw)<l(w)$ then by Lemma~\ref{nc lemma 2},~$s$ is a cover reflection of~$w$.
In this case, $s\in sS's$ by Lemma~\ref{s in}, so that $s\in S'$, contradicting the fact that $\nu_c(W')>0$.
This contradiction shows that $l(sw)>l(w)$ so that $sw$ is $c$-sortable by Lemma~\ref{scs}.
If~$s$ is a cover reflection of $sw$ then by Lemma~\ref{nc lemma} the remaining cover reflections of $sw$ are in~$W_{\br{s}}$.
Thus the cover reflections of~$w$ are all in $sW_{\br{s}}s$, so that \mbox{$y=\nc_{scs}(w)\in sW_{\br{s}}s$}.
This contradicts our assumption that $x=sys$ is not in~$W_{\br{s}}$ for any initial letter~$s$ of~$c$.
Thus~$s$ is not a cover reflection of $sw$, so $\nc_c(sw)=s\cdot\nc_{scs}(w)\cdot s=x$.

If~$w'$ is any other $c$-sortable element with $\nc_c(w')=x$ then by Lemmas~\ref{cover para} and~\ref{abs para}, $w'\not\in W_{\br{s}}$.
Thus $l(sw')<l(w')$ by Lemma~\ref{sc}.
Also,~$s$ is not a cover reflection of~$w'$.
If it were then~$s$ would be in~$W'$ and therefore by Lemma~\ref{s in} in~$S'$, contradicting the fact that $\nu_c(W')>0$.
Thus $\nc_{scs}(sw')=s\cdot\nc_c(w')\cdot s=y$ so that $sw'=w$ and therefore $w'=sw$.
\end{proof}

\begin{remark}\rm
%sinceV1:  fixed typo ``biproduct''
As a byproduct of Theorem~\ref{nc}, we obtain a canonical reduced $T$-word for every~$x$ in $[1,c]_T$.
The letters are the canonical generators of the associated parabolic subgroup, or equivalently the cover reflections of $(\nc_c)^{-1}(x)$,
occurring in the order induced by the $c$-sorting word for $(\nc_c)^{-1}(x)$.
This is canonical, up to the choice of reduced word for~$c$.
Changing the reduced word for~$c$ alters the choice of canonical reduced $T$-word for~$x$ by commutations of letters.

%sinceV1:  made a more definite statement about the connection with ABW:
In~\cite{ABW}, it is shown that for~$c$ bipartite, the natural labeling of $[1,c]_T$ is an EL-shelling with respect to the reflection order obtained from what we here call the $c$-sorting word for~$w_0$.
In particular, the labels on the unique maximal chain in $[1,x]_T$ constitute a canonical reduced $T$-word for~$x$.
%rev:  Check exact reference to ABW result in following sentence:
It is apparent from the proof of \cite[Theorem 3.5(ii)]{ABW} that these two choices of canonical reduced $T$-words are identical in the bipartite case, for $W$ crystallographic.
Presumably the same is true for non-crystallographic $W$.
\end{remark}

%sinceV1:  added the following remark to aid in writing sort_camb:
\begin{remark}\rm
Theorem~\ref{nc} can also be rephrased in terms of the set $\cov(w)$ of cover reflections of a $c$-sortable element $w$.
Specifically, the map $\cov$ takes a $c$-sortable element to the set of canonical generators of the parabolic subgroup associated to some noncrossing partition (with respect to $c$).
The content of Theorem~\ref{nc} is that this map is a bijection:  For any parabolic subgroup $W'$ associated to a noncrossing partition (with respect to $c$), there is a unique $c$-sortable element $w$ such that $\cov(w)$ is the set of canonical generators of $W'\!$.\@

Also, in light of Theorem~\ref{nc}, basic facts about noncrossing partitions have analogous statements in terms of sets of cover reflections of $c$-sortable elements.
For example, the fact that (i) implies (iv) in Lemma~\ref{abs para} becomes the following statement:
If~$s$ is initial in~$c$ and $w\in W_{\br{s}}$ is $sc$-sortable then there exists a $c$-sortable element $w'$ with $\cov(w')=\set{s}\cup\cov(w)$.
To see this, let $x=\nc_c(w)$ and let $W'$ be the parabolic subgroup associated to $x$.
Then $sx$ is a noncrossing partition by Lemma~\ref{abs para}.
The canonical generators of the parabolic subgroup associated to $sx$ are $\set{s}$ union the canonical generators of $W'$.
Thus $w'=\nc^{-1}(sx)$ is the desired $c$-sortable element.
\end{remark}

\section{Clusters}
\label{cl prelim}
In this section we define $c$-clusters, a slight generalization (from crystallographic Coxeter groups to all finite Coxeter groups) of the $\Gamma$-clusters of~\cite{MRZ}.
These in turn generalize the clusters of~\cite{ga}.
We build the theory of clusters within the framework of Coxeter groups, rather than in the framework of root systems.
This is done in order to avoid countless explicit references to the map between positive roots and reflections in what follows. 
Readers familiar with root systems will easily make the translation to the language of {\em almost positive roots} of~\cite{ga} 
and~\cite{MRZ}.

Let $-S$ denote the set $\set{-s:s\in S}$ of formal negatives of the simple generators of~$W,$\ and let $\Tge$ be $T\cup(-S)$.
(Recall that $T$ is the set of all reflections of~$W.$)
For each $J\subseteq S$, the notation $T_J$ stands for $T\cap W_J$ and $(T_J)_{\ge-1}$ denotes $T_J\cup(-J)$.

For each $s\in S$, define an involution $\sigma_s:\Tge\to\Tge$ by
\[\sigma_s(t):=\left\lbrace\begin{array}{ll}
-t&\mbox{if }t=\pm s,\\
t&\mbox{if }t\in (-S)\mbox{ and }t\neq -s,\mbox{ or}\\
sts&\mbox{if }t\in T-\set{s}.
\end{array}\right.\]

Given any $t\in\Tge$, any Coxeter element~$c$ and any sequence $s_0,s_1,\ldots,s_k$ of simple generators, define $t_0=t$, $c_0=c$,
$t_i=\sigma_{s_{i-1}}(t_{i-1})$ and $c_i=s_{i-1}c_{i-1}s_{i-1}$.
\begin{lemma}
\label{mu lemma}
For any $t\in\Tge$ there exists a sequence $s_0,s_1,\ldots,s_k$ with $s_i$ initial in $c_i$ for $i=0,1,\ldots,k$ and $t_k=-s_k$.
\end{lemma}
\begin{proof}
If $t=-s$ for some initial letter~$s$ of~$c$ then the desired sequence has a single entry~$s$.
If $t=-s$ for some $s\in S$ which is not initial in~$c$ then let $c'$ be any Coxeter element for which~$s$ is initial.
By Lemma~\ref{swap sources}, $c'$ can be obtained from~$c$ by a sequence of conjugations, each by an initial letter of the current Coxeter
element.
Thus some sequence $s_0,\ldots,s_k$ has $t_k=t_{k-1}=\cdots=t_0=-s_k=-s$ with $s_i$ initial in $c_i$ for $i=0,1,\ldots,k$.

If $t\not\in-S$ then let~$W'$ be the two-element parabolic subgroup generated by $t$, so that the set of canonical generators of~$W'$ is 
$S'=\set{t}$.
Lemma~\ref{nu lemma} guarantees the existence of a sequence $s_0,s_1,\ldots,s_j$ with $s_i$ initial in $c_i$ for $i=0,1,\ldots,k$ 
and $t_i=s_{i-1}t_{i-1}s_{i-1}=\sigma_{s_{i-1}}(t_{i-1})$ for each $i\in[k]$, such that $t_k=s_k$.
Then $t_{k+1}=-s_k$.
We append a sequence of letters, as in the previous paragraph, which make $s_k$ initial while preserving $-s_k$, then append the letter $s_k$ as the final letter of the sequence.
\end{proof}
Define $\mu_c(t)$ to be the smallest $k$ such that there exists a sequence $s_0,s_1,\ldots,s_k$ satisfying the conclusions of 
Lemma~\ref{mu lemma}.
In particular $\mu_c(t)=0$ if and only if $t=-s$ for an initial letter~$s$ of~$c$.

We now define a symmetric binary relation $\cm_c$ called the {\em $c$-compatibility} relation.
\begin{prop}
\label{compat}
There exists a unique family of symmetric binary relations $\cm_c$ on $\Tge$, indexed by Coxeter elements~$c$, with the 
following properties:
\begin{enumerate}
\item[(i) ]For any $s\in S$, $t\in\Tge$ and Coxeter element~$c$,
\[-s\cm_ct\mbox{ if and only if }t\in(T_{\br{s}})_{\ge-1}.\]
\item[(ii) ]For any $t_1,t_2\in\Tge$ and any initial letter~$s$ of~$c$, 
\[t_1\cm_c t_2\mbox{ if and only if }\sigma_s(t_1)\,\cm_{scs}\,\sigma_s(t_2).\]
\end{enumerate}
\end{prop}
By symmetry, (ii) is equivalent to the same statement with a final letter~$s$ of~$c$ in place of an initial letter.
\begin{proof}
In light of Lemma~\ref{mu lemma}, such a family is unique if it exists.
It is sufficient to verify existence in the case where~$W$ is irreducible.
Marsh, Reineke and Zelevinsky construct a family of relations for Coxeter groups of types A, D and E satisfying
(i) and (ii) \cite[Propositions 3.4 and 4.10]{MRZ}.
They also point out (see \cite[Section~1]{MRZ}) that this construction and properties can be extended to all finite crystallographic Coxeter groups by standard ``folding'' techniques.

Verification of existence for the remaining types does not rely on any results of Section~\ref{cl}, but is informed by them, as we now explain.
Theorem~\ref{cl} implies that for $t_1,t_2\in \Tge$, $t_1\cm_c t_2$ if and only if there exists a $c$-sortable element~$w$ such that $t_1,t_2\in\cl_c(w)$, where $\cl_c(w)$ is defined in Section~\ref{cl sec} using only results of Section~\ref{sort sec} and earlier.
Taking this as a definition of $c$-compatibility, one easily verifies (i) and (ii) for Coxeter groups of rank two, and checks (i) and (ii) by computer for~$H_3$ and~$H_4$.
%To do this in sort_maple, run:
%crys:= false: for R in [H3,H4] do check_compat_all(R) od;
\end{proof}

\begin{prop}
\label{restrict}
Let $J\subseteq S$ and $t_1,t_2\in(T_J)_{\ge-1}$ and let $c'$ be the Coxeter element for~$W_{\br{s}}$ obtained from~$c$ by restriction.
Then $t_1\cm_ct_2$ if and only if $t_1\cm_{c'}t_2$.
\end{prop}
\begin{proof}
This is \cite[Proposition~3.4]{MRZ} in the crystallographic case, so it suffices to check the non-crystallographic irreducible cases.
The rank-two case is trivial and $H_3$ and $H_4$ are checked by computer. 
\end{proof}
%To do this in sort_maple, run:
%crys:= false: for R in [H3,H4] do check_compat_res_all(R) od;

\begin{prop}
\label{c inv}
The relations $\cm_c$ and $\cm_{c^{-1}}$ coincide.
%For any Coxeter element~$c$, $t_1\cm_ct_2$ if and only if $t_1\cm_{c^{-1}}t_2$.
\end{prop}
\begin{proof}
Induction on $\mu_c$, using properties (i) and (ii) of $c$-compatibility.
\end{proof}
The crystallographic case of Proposition~\ref{c inv} is \cite[Proposition~3.1]{MRZ}.
In \cite[Corollary 4.12]{MRZ}, it is pointed out that for a bipartite Coxeter element~$c$, the relation $\cm_c$ coincides with the notion of compatibility defined in~\cite{ga}.
This is also seen by comparing properties (i) and (ii) of $c$-compatibility and Proposition~\ref{c inv} with the construction of \cite[Section~3]{ga}.

A {\em $c$-compatible subset} of $\Tge$ is a set of pairwise $c$-compatible elements of $\Tge$.
A {\em $c$-cluster} is a maximal $c$-compatible subset.
A $c$-cluster is called {\em positive} if it contains no element of $-S$.

\section{Sortable elements and clusters}
\label{cl sec}
The main result of this section is a bijection between $c$-sortable elements and $c$-clusters.
Let~$w$ be a $c$-sortable element with $c$-sorting word $a_1a_2\cdots a_k$.
If $s\in S$ occurs in $a$ then the \emph{last reflection} for~$s$ in~$w$ is $a_1a_2\cdots a_ja_{j-1}\cdots a_2a_1$, where $a_j$ is the rightmost occurrence of~$s$ in $a$.
If~$s$ does not occur in $a$ then the last reflection for~$s$ in~$w$ is the formal negative $-s$.
Let $\cl_c(w)$ be the set of last reflections of~$w$.
This is an $n$-element subset of $\Tge$.
This map does not depend on the choice of reduced word for~$c$, because any two $c$-sorting words for~$w$ are related by commutations of simple generators.

\begin{theorem}
\label{cl}
The map $w\mapsto\cl_c(w)$ is a bijection from the set of $c$-sortable elements to the set of $c$-clusters.
Furthermore, $\cl_c$ restricts to a bijection between $c$-sortable elements with full support and positive $c$-clusters.
\end{theorem}

\begin{example}\rm
\label{B2cl}
We continue the example of $W=B_2$ and $c=s_0s_1$.
Clusters in $B_2$ correspond to collections of diagonals which define centrally symmetric triangulations of the hexagon shown below.
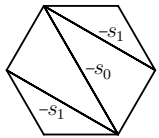
\[\scalebox{.7}{\epsfbox{B2seed.ps}}\]
Each element of $\Tge$ is represented by a diameter or a centrally symmetric pair of diagonals.
For details, see \cite[Section~3.5]{ga}.
Figure~\ref{B2clfig} illustrates the map $\cl_c$ on $c$-sortable elements for this choice of~$W$ and~$c$.
\end{example}

{\small
\begin{figure}[ht]
\[\begin{array}{|c||c|c|c|c|c|c|}\hline
&&&&&&\\[-7 pt]
w	&1		&s_0		&s_0s_1		&s_0s_1|s_0	&s_0s_1|s_0s_1	&s_1\\[3 pt]\hline
&&&&&&\\[-7 pt]
\cl_c(w)	&-s_0, -s_1	&s_0, -s_1	&s_0, s_0s_1s_0	&s_1s_0s_1, s_0s_1s_0&s_1s_0s_1, s_1	&-s_0, s_1\\[3pt]\hline
&&&&&&\\[-7 pt]
	&\epsfbox{B2cl.1.ps}&\epsfbox{B2cl.a.ps}	&\epsfbox{B2cl.ab.ps}&\epsfbox{B2cl.aba.ps}&\epsfbox{B2cl.abab.ps}&\epsfbox{B2cl.b.ps}\\
\hline
\end{array}\]
\caption{The map $\cl_c$}
\label{B2clfig}
\end{figure}
}

\begin{example}\rm
By way of contrast with Example~\ref{Annc}, we offer no characterization of the map $\cl_c$ on permutations satisfying (A), even in the $231$-avoiding case.
Such~a characterization is not immediately apparent, due to the dependence of $\cl_c(w)$ on~a specific choice of reduced word for $w$.
\end{example}

\begin{remark}\rm
\label{inf cl}
Even for infinite $W,$ the map $\cl_c$ associates to each $c$-sortable element an $n$-element subset of $\Tge$.
However, for infinite $W,$ it is not even clear how $c$-compatibility should be defined, and in particular the proofs in this section apply to the finite case only. 
As mentioned in the proof of Proposition~\ref{compat}, Theorem~\ref{cl} implies the following characterization of $c$-compatibility:
Distinct elements $t_1$ and $t_2$ of $\Tge$ are $c$-compatible if and only if there exists a $c$-sortable element~$w$ such that $t_1,t_2\in\cl_c(w)$.
Thus the map $\cl_c$ itself might conceivably provide some insight into compatibility in the infinite case.
\end{remark}

The strategy for proving Theorem~\ref{cl} is the same as for Theorem~\ref{nc}, but with fewer complications.
We argue by induction on rank and length using Lemmas~\ref{sc} and~\ref{scs} and the following two lemmas.

\begin{lemma}
\label{cl s}
Let~$s$ be initial in~$c$ and let~$w$ be $c$-sortable.
If $l(sw)>l(w)$ then $\cl_c(w)=\set{-s}\cup\cl_{sc}(w)$.
If $l(sw)<l(w)$ then $\cl_c(w)=\sigma_s(\cl_{scs}(sw))$.
\end{lemma}
\begin{proof}
Use a reduced word for $c$ whose first letter is $s$ to write a $c$-sorting word for $w$.
Delete the first letter from this word to obtain an $scs$-sorting word for $sw$.
Comparing these two words, the lemma is immediate.
\end{proof}

\begin{lemma}
\label{cl lemma}
Let~$s$ be final in~$c$ and let~$w$ be $c$-sortable.
Then $l(sw)<l(w)$ if and only if $s\in\cl_c(w)$.
\end{lemma}
\begin{proof}
The ``if'' direction holds by definition.
Suppose~$s$ is final in~$c$ and $l(sw)<l(w)$ for a $c$-sortable element~$w$ with $c$-sorting word $a_1a_2\cdots a_k$.
Let $i$ be such that $s=a_1a_2\cdots a_ia_{i-1}\cdots a_1$.
We now prove, by induction on $k-i$, the three claims below.
Note that claim (ii) is the lemma.
\begin{enumerate}
\item[(i) ]$sw=wa_i$, 
\item[(ii) ]$a_j\neq a_i$ for $i<j\le k$, and
\item[(iii) ]$a_j$ commutes with $a_i$ for $i<j\le k$.
\end{enumerate}
If $i=k$ then claim (i) holds by construction and claims (ii) and (iii) are vacuous.
Suppose $i<k$ and let $w'=wa_k$.
Then~$w'$ is $c$-sortable with $c$-sorting word $a_1a_2\cdots a_{k-1}$.
By induction the three claims hold for~$w',$\ so that in particular $a_j\neq a_i$ for $i<j<k$.
If $a_k=a_i$ then since $a_i$ commutes with $a_j$ for $i<j< k$, the word $a_1a_2\cdots a_k$ is not reduced.
This contradiction shows that $a_k\neq a_i$ so claim (ii) is established.
Claim (i) for~$w'$ says that $sw'=w'a_i$ so that $sw=w(a_ka_ia_k)$.
By Lemma~\ref{nc lemma 2}, $sw=ws'$ for some simple generator $s'$ and therefore $a_i$ must commute with $a_k$, so that $s'=a_i$.
This establishes claims (i) and (iii) for~$w$.
\end{proof}

We now prove the main theorem of the section.

\begin{proof}[Proof of Theorem~\ref{cl}]
We first verify that for any $c$-sortable element~$w$, the set $\cl_c(w)$ is a $c$-cluster.
Let~$w$ be a $c$-sortable element and let~$s$ be an initial letter of~$c$.
If $l(sw)>l(w)$ then~$w$ is an $sc$-sortable element of~$W_{\br{s}}$ by Lemma~\ref{sc}.
By induction on the rank of~$W,$\ $\cl_{sc}(w)$ is an $sc$-cluster.
Proposition~\ref{restrict} implies that $\cl_{sc}(w)$ is a set of pairwise $c$-compatible elements of $\Tge$.
This set is maximal with respect to pairwise $c$-compatibility and inclusion in $(T_J)_{\ge-1}$.
By Lemma~\ref{cl s} and property (i) of $c$-compatibility, $\cl_c(w)$ is a $c$-cluster.
If $l(sw)<l(w)$ then by Lemma~\ref{cl s}, $\cl_c(w)=\sigma_{s}(\cl_{scs}(sw))$.
By induction on the length of~$w$, $\cl_{scs}(sw)$ is an $scs$-cluster, so by property (ii) of $c$-compatibility, $\cl_c(w)$ is a $c$-cluster.

Let~$C$ be a cluster and let $\mu_c(C)$ be the minimum of $\mu_c(t)$ for $t\in C$.
We now prove by induction on $\mu_c(C)$ that there exists a unique $c$-sortable $w\in W$ with $\cl_c(w)=C$.
If $\mu_c(C)=0$ then there is some $-s\in C$ such that~$s$ is initial in~$c$.
In particular any $c$-sortable element~$w$ with $\cl_c(w)=C$ must have $w\in W_{\br{s}}$.
By Lemma~\ref{cl s}, $\cl_{sc}(w)=C-\set{-s}$, which by Proposition~\ref{restrict} is an $sc$-cluster in~$W_{\br{s}}$.
By induction on the rank of~$W,$\ there exists a unique $sc$-sortable element $w\in W_{\br{s}}$ such that $\cl_{sc}(w)=C-\set{-s}$.
This element is $c$-sortable by Lemma~\ref{sc} and has $\cl_c(w)=C$ by Lemma~\ref{cl s}.

If $\mu_c(C)>0$ then let $s_0,s_1,\ldots,s_{\mu_c(C)}$ be a sequence satisfying the conditions of Lemma~\ref{mu lemma} for some $t\in C$ with $\mu_c(t)=\mu_c(C)$.
Let $s=s_0$ and consider the $scs$-cluster $C'=\sigma_{s}(C)$.
Then $\mu_{scs}(C')=\mu_c(C)-1$, so by induction, there is a unique $scs$-sortable element~$w$ such that $\cl_{scs}(w)=C'$.
If $l(sw)<l(w)$ then by Lemma~\ref{cl lemma}, $s\in\cl_{scs}(w)=C'$.
But this means that $\sigma_s(s)\in\sigma_s(C')$, or in other words, $-s\in C$, contradicting the assumption that $\mu_c(C)>0$.
By this contradiction we conclude that $l(sw)>l(w)$ and therefore by Lemma~\ref{scs}, $sw$ is a $c$-sortable element.
Lemma~\ref{cl s} says that $\cl_c(sw)=\sigma_s(\cl_{scs}(sw))=C$.
Let~$x$ be any $c$-sortable element with $\cl_c(x)=C$.
Since $-s\not\in C$ we have $l(sx)<l(x)$ so that $sx$ is an $scs$-sortable element with $\cl_{scs}(sx)=\sigma_s(C)=C'$.
By the uniqueness of~$w$, $sx=w$, so that $x=sw$.
Thus $sw$ is the unique $c$-sortable element mapping to~$C$ under the map $\cl_c$.

The last statement of the theorem is immediate from the definition of $\cl_c$.
\end{proof}

\section{Enumeration}
\label{enum}
In this section we briefly discuss the enumeration of sortable elements.
The {\em~$W$-Catalan number} is given by the following formula, in which $h$ is the Coxeter number of~$W$ and the $e_i$ are the exponents of~$W.$\
\[\Cat(W)=\prod_{i=1}^n \frac{e_i + h + 1}{e_i + 1} \,. \]
The values of the~$W$-Catalan number for finite irreducible $W$ are tabulated below.\\[6 pt]
\noindent
\scalebox{.95}{
\mbox{
$
\begin{array}{|c|c|c|c|c|c|c|c|c|c|c|c|}
\hline
A_n & B_n & D_n & E_6 & E_7 & E_8 & F_4  & G_2 & H_3 & H_4 & I_2(m)\\\hline
&&&&&&&&&&\\[-.1in]
    \textstyle\frac{1}{n+2} \binom{2n+2}{n+1}
& \textstyle\binom{2n}{n} 
& \textstyle\frac{3n-2}{n}\binom{2n-2}{n-1} 
& \textstyle 833
& \textstyle 4160
& \textstyle 25080
& \textstyle 105
& \textstyle 8
& \textstyle 32
& \textstyle 280
& \textstyle m+2 \\[.05in]
\hline
\end{array}
$
}
}\\[2 pt]

The noncrossing partitions (with respect to any Coxeter element) in an irreducible finite Coxeter group~$W$ are counted by the $W$-Catalan number \cite{Bessis,picantin,Reiner}.
The $c$-clusters (for any Coxeter element~$c$) of an irreducible finite Coxeter group~$W$ are also counted by $\Cat(W)$.
This follows from \cite[Corollary 4.11]{MRZ} and \cite[Proposition~3.8]{ga} for the crystallographic case, or is proved in any finite case by combining Theorems~\ref{nc} and~\ref{cl}.
We refer the reader to \cite[Section~5.1]{rsga} for a brief account of other objects counted by the~$W$-Catalan number.
By Theorem~\ref{nc} or Theorem~\ref{cl} we have the following.

\begin{theorem}
\label{num sort}
For any Coxeter element~$c$ of~$W,$\ the $c$-sortable elements of~$W$ are counted by the~$W$-Catalan number.
\end{theorem}

The {\em positive~$W$-Catalan number} is the number of positive $c$-clusters ($c$-clusters containing no element of $-S$).
The following is an immediate corollary of Theorem~\ref{cl}.

\begin{cor}
\label{positive}
For any Coxeter element~$c$, the $c$-sortable elements not contained in any proper standard parabolic subgroup are counted by the positive~$W$-Catalan number. 
\end{cor}
The map $\nc_c$ also respects this positive $W$-Catalan enumeration:
By Lemmas~\ref{cover para} and~\ref{abs para}, the map $\nc_c$ maps the sortable elements not contained in any proper standard parabolic subgroup to the noncrossing partitions not contained in any proper standard parabolic subgroup.

The~$W$-Narayana numbers count noncrossing partitions by their rank.
That is, the $k$th~$W$-Narayana number is the number of elements of $[1,c]_T$ whose absolute length is $k$.
The following is an immediate corollary of Theorem~\ref{nc}.

\begin{cor}
\label{Narayana}
For any Coxeter element~$c$, the $c$-sortable elements of~$W$ which have exactly $k$ descents are counted by the $k$th~$W$-Narayana number.
\end{cor}

\begin{remark}\rm
\label{cl Narayana}
The $k$th~$W$-Narayana number is also the $k$th component in the $h$-vector of the simplicial~$W$-associahedron.
Using results from~\cite{con_app} and~\cite{sort_camb}, one associates a complete fan to $c$-sortable elements.
This fan has the property that any linear extension of the weak order on $c$-sortable elements is a shelling.
%sinceV1:  The following sentences altered.
In~\cite{camb_fan}, David Speyer and the author show that the map $\cl_c$ induces a combinatorial isomorphism.
Thus as a special case of a general fact explained in the discussion following \cite[Proposition~3.5]{con_app}, the $h$-vector of $\Delta_c$ has entry $h_k$ equal to the number of $c$-sortable elements with exactly $k$ descents.
This gives an alternate proof of Corollary~\ref{Narayana} and, by composing bijections, a bijective explanation of why counting noncrossing partitions by rank recovers the $h$-vector of the~$W$-associahedron.
\end{remark}

\section*{Acknowledgments}
I am grateful to Sara Billey, Sergey Fomin, Jon McCammond, Vic Reiner, John Stembridge, 
Hugh Thomas and Andrei Zelevinsky for helpful conversations.
Thanks also to Stembridge for advice on efficient computations. 

The computations described here were done with \texttt{Maple} using John Stembridge's \texttt{coxeter/weyl} packages.
The \texttt{Maple} source code is currently available at\\
\texttt{http://www.math.lsa.umich.edu/\textasciitilde nreading/papers/sortable/}.\\

\newcommand{\journalname}[1]{\textrm{#1}}
\newcommand{\booktitle}[1]{\textrm{#1}}

\end{document}